\documentclass[reqno]{amsart}
\pdfoutput=1

\usepackage{amsmath,amsthm,amssymb,amsfonts,amscd,amsxtra}
\usepackage[]{hyperref}
\usepackage{graphicx}
\usepackage{mathrsfs} 
\usepackage[all]{xy}  
\usepackage{stmaryrd} 
\usepackage{latexsym} 
\usepackage{caption}

\newtheorem{thm}{Theorem}[section]
\newtheorem{prop}[thm]{Proposition}
\newtheorem{lem}[thm]{Lemma}
\newtheorem{cor}[thm]{Corollary}

\newtheorem*{claim*}{Claim}

\newtheorem{step}{Step}[]

\theoremstyle{definition}
\newtheorem{defn}{Definition}[section]
\newtheorem{example}{Example}

\newtheorem*{examples*}{Examples}
\newtheorem*{example*}{Example}

\theoremstyle{remark}
\newtheorem{rmk}[thm]{Remark}

\newcommand{\deco}{{}}
\newcommand{\sesLS}{\sharp}
\newcommand{\Fuk}{\mathscr F}
\newcommand{\FukExt}{\mathscr F^\sesLS}
\newcommand{\cobGr}{\Omega_{\mathrm{Lag}}}
\newcommand{\cobGrExt}{\Omega_{\mathrm{Lag}}^\sesLS}
\newcommand{\Lags}{\mathscr L}
\newcommand{\LagsExt}{\mathscr L^\sesLS}

\title{The Lagrangian cobordism group of $T^2$}
\author{Luis Haug}
\address{Department of Mathematics, ETH Z\"urich, R\"amistrasse 101,
    8092 Z\"urich, Switzerland}
\email{haug@math.ethz.ch}
\date{\today}


\begin{document}
\begin{abstract}
    We compute the Lagrangian cobordism group of the standard
    symplectic 2-torus and show that it is isomorphic to the
    Grothendieck group of its derived Fukaya category. The proofs use
    homological mirror symmetry for the 2-torus.
\end{abstract}
\maketitle
\tableofcontents

\section{Introduction}
\label{sec:introduction}

\begin{figure}[t]
    \centering
    \includegraphics[scale=1.3]{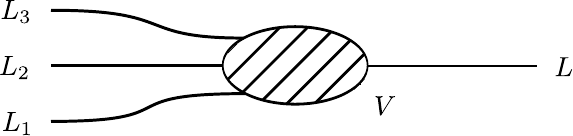}
    \caption{Projection of a Lagrangian cobordism $V: L \leadsto
        (L_1,L_2,L_3)$.}
    \label{fig:cobordism}
\end{figure}

The Fukaya category $\Fuk(M)$ of a symplectic manifold $(M,\omega)$ is
an $A_\infty$-category whose objects are the Lagrangian submanifolds
of $M$, and whose morphism spaces are Floer cochain groups. Its
derived category $D \Fuk(M)$ is triangulated, and thus it possesses a
mechanism for generating objects by taking cones of morphisms. Recent
work of Biran--Cornea
\cite{Biran-Cornea--Lag-Cob-I,Biran-Cornea--Lag-Cob-II} provides a way
of understanding cone decompositions of objects in $D \Fuk(M)$
geometrically via Lagrangian cobordisms. According to their
definition, a Lagrangian cobordism $V:(L_1,\dots,L_r) \leadsto
(L_1',\dots,L_s')$ between tuples of Lagrangians in $M$ is a
Lagrangian submanifold $V$ of $\mathbb R^2 \times M$ with cylindrical
ends corresponding to the $L_i$ and $L_j'$. Figure \ref{fig:cobordism}
displays the projection of such a cobordism to $\mathbb R^2$. The main
result of \cite{Biran-Cornea--Lag-Cob-I,Biran-Cornea--Lag-Cob-II} is
that a Lagrangian cobordism of the form
\begin{equation*}
    V: L \leadsto (L_1,\dots,L_s)
\end{equation*}
leads to an iterated cone decomposition of $L$ in $D \Fuk(M)$ whose
``building blocks'' are the $L_i$.

Partial information about the triangulated structure of $D \Fuk(M)$ is
captured by its \emph{Grothendieck group} $K_0(D \Fuk(M))$. It is
generated by the objects of $\Fuk(M)$, that is, the Lagrangians in
$M$, with relations coming from exact triangles in $D \Fuk(M)$. On the
cobordism side, one can naturally define a \emph{Lagrangian cobordism
    group} $\cobGr(M)$. It is generated by the Lagrangians in $M$, or
a suitable subset thereof, and has relations coming from Lagrangian
cobordisms. As an immediate consequence of Biran--Cornea's results,
there exists a surjective group homomorphism
\begin{equation*}
    \Theta: \cobGr(M) \to K_0(D \Fuk(M)),
\end{equation*}
induced by $L \mapsto L$. It is natural to ask if $\Theta$ is an
isomorphism, and if not, what its kernel is. This question formalises
the question to what extent the triangulated structure of $D \mathscr
F(M)$ can be explained geometrically by the existence of Lagrangian
cobordisms. In general, nothing is known about the kernel of $\Theta$.

\subsection*{Main results.}
\label{sec:main-result}
The case we consider is when $(M,\omega)$ is the standard symplectic
2-torus $(T^2 = \mathbb R^2/\mathbb Z^2, \omega_{\mathrm{std}} = dx
\wedge dy)$. The version of the Lagrangian cobordism group
$\cobGr(T^2)$ we study is generated by non-contractible simple closed
curves in $T^2$. 

Our first main result is the following statement.
\begin{thm}
    \label{thm:MainThm}
    The natural group homomorphism
    \begin{equation*}
        \Theta: \cobGr(T^2) \to K_0(D \mathscr F(T^2))
    \end{equation*}
    is an isomorphism.
\end{thm}

Theorem \ref{thm:MainThm} has implications of three different
kinds. First, as already indicated above, it tells us that the set of
relations in $K_0(\Fuk(T^2))$ is generated by ones coming from
Lagrangian cobordisms. Second, together with Theorem
\ref{thm:MainThm-SES} it can be regarded as a computation of
$K_0(D\Fuk(T^2))$, which to the author's knowledge has not been
carried out before (but cf. \cite{Abouzaid--Fuk-higher-genus-surfaces}
for a computation of $K_0(D\Fuk(\Sigma))$ for higher genus
surfaces). Third, it gives information about Lagrangian cobordisms in
$\mathbb R^2 \times T^2$. There are two known constructions of such
cobordisms based on Hamiltonian isotopy and on Lagrangian surgery,
but we do not know if the resulting cobordisms are the only ones that
exist. Theorem \ref{thm:MainThm} does not rule out that there are
more, but it (or rather its proof) shows that the known cobordisms
generate the set of relations in $\cobGr(T^2)$.

The second main result answers the question what $\cobGr(T^2)$ looks
like:
\begin{thm}
    \label{thm:MainThm-SES}
    There exists a canonical short exact sequence
    \begin{equation*}
        0 \to \mathbb R/\mathbb Z \xrightarrow{\zeta} \cobGr(T^2) \xrightarrow{\eta}
        H_1(T^2;\mathbb Z)\to 0.
    \end{equation*}
\end{thm}

The map $\eta$ is the obvious one, given by $[L]_{\Omega} \mapsto
[L]_{H_1}$. The map $\zeta$ takes $x \in \mathbb R/\mathbb Z$ to the
class represented by the boundary of a cylinder of area $x$. In other
words, $\zeta(x)$ is represented by the difference $L - L'$, with $L$
any non-contractible curve, and $L'$ any isotopic curve such that the
area swept out during any isotopy from $L'$ to $L$ is $x$ (this area
is determined by the two curves up to an integer). Since two such
curves are Hamiltonian isotopic if and only if they bound a cylinder
of area 0, the injectivity of $\zeta$ says that $\cobGr(T^2)$
distinguishes different Hamiltonian isotopy classes.

While most of the proof of Theorem \ref{thm:MainThm-SES} is quite
elementary, it is not clear how to rule out in a direct way possible
cobordisms that might obstruct the injectivity of $\zeta$. Our proof
relies on the connection to homological mirror symmetry discussed
below.

\subsection*{Proving Theorem \ref{thm:MainThm} using HMS}
\label{sec:towards-proof-theor}
We will first prove an extended version of Theorem \ref{thm:MainThm}
for which we consider a Fukaya category $\FukExt(T^2)$ defined over a
Novikov field $\Lambda$ (consisting of formal power series with
$\mathbb C$-coefficients), whose objects are Lagrangians in $T^2$ that
are decorated with certain local systems of $\Lambda$-vector
spaces.\footnote{We remark that from certain perspectives
    $\FukExt(T^2)$ might actually be a more natural version of the
    Fukaya category to look at than $\Fuk(T^2)$.} The definition of
the cobordism group is also modified accordingly: The Lagrangians as
well as the cobordisms carry local systems, and there are additional
relations coming from short exact sequences of local systems. We
denote the resulting group by $\cobGrExt(T^2)$. As in the case without
local systems, there is a natural surjective group homomorphism
$\Theta^\sesLS: \cobGrExt(T^2) \to K_0(D \FukExt(T^2))$.

\begin{thm}
    \label{thm:MainThmExtended}
    The natural group homomorphism
    \begin{equation*}
        \Theta^\sesLS: \cobGrExt(T^2) \to    K_0(D \FukExt(T^2))
    \end{equation*}
    is an isomorphism. 
\end{thm}

As the main ingredient in the proof of Theorem
\ref{thm:MainThmExtended}, apart from Biran--Cornea's theory, we use
that $T^2$ is one of the symplectic manifolds for which the
homological mirror symmetry conjecture has been proven. The statement
of relevance to us is the result of Abouzaid--Smith
\cite{Abouzaid-Smith--HMS4torus}, who on their way to HMS for $T^4$
construct a triangulated equivalence $D^b(X) \simeq D^\pi \Fuk(T^2)$
between the derived category of coherent sheaves of an elliptic curve
$X$ defined over the Novikov field $\Lambda$ and the split-closed
derived Fukaya category of $T^2$. (A more refined version of HMS for
$T^2$ was recently proven by Lekili--Perutz
\cite{Lekili-Perutz--Arithmetic-MS-T2}, but this is not needed for our
purposes.)  

One can adapt Abouzaid--Smith's result to our setting, such as to
obtain a triangulated equivalence
\begin{equation*}
    D^b(X) \simeq D \FukExt(T^2).
\end{equation*}
(In particular, this shows that taking the split-closure of $D
\Fuk(T^2)$ or using appropriate local systems is equivalent here.)
What needs to be checked is that Abouzaid--Smith's functor, which is
given explicitly only on a small collection of split-generators, takes
every sheaf to a Lagrangian with a local system (as opposed to some
summand of some non-trivial complex). We do so by matching up sheaves
and Lagrangians with local systems in an inductive manner, using
Atiyah's classification of vector bundles on elliptic curves
\cite{Atiyah--VecBunEllCurve--1957} respectively surgery and
cobordisms to understand the structures of $D^bX$ and
$D\FukExt(T^2)$. The result of this procedure is stated in Proposition
\ref{prop:recov-mirr-funct}.

The resulting isomorphism between Grothendieck groups
\begin{equation*}
    K_0(D^b(X)) \cong K_0(D \FukExt(T^2))
\end{equation*}
allows us to understand relations in $K_0(D \FukExt(T^2))$ via the
well-understood group $K_0(D^b(X))$, and to check that the ``obvious''
inverse to the map $\Theta^\sesLS: \cobGrExt(T^2) \to K_0(D
\FukExt(T^2))$ is well-defined. This proves Theorem
\ref{thm:MainThmExtended}, from which Theorem \ref{thm:MainThm}
follows.

\subsection*{Organisation of the paper.}
\label{sec:organisation-paper}
In Section \ref{sec:lagr-cobord-group}, we recall the definition of
Lagrangian cobordisms and define the groups $\cobGr(M)$ and
$\cobGrExt(M)$, and Section \ref{sec:Fuk-T2} serves to describe the
Fukaya category $\FukExt(T^2)$. In Section \ref{sec:lagr-cobord-cones}
we explain how Lagrangian cobordisms lead to iterated cone
decompositions in the derived Fukaya category, focusing on the (small)
modifications necessary to make Biran--Cornea's proofs work in our
setting. In Section \ref{sec:lagr-cobord-T2}, we describe in detail
the Lagrangian cobordisms resulting from surgering curves in surfaces
and discuss some specifics in the case of $T^2$. Section
\ref{sec:comp-cobGr} serves to explain and prove as much as we can
about $\cobGr(T^2)$ without any mirror symmetry considerations; in
particular, we prove Theorem \ref{thm:MainThm-SES} modulo the
injectivity of the map $\zeta: \mathbb R/\mathbb Z \to
\cobGr(T^2)$. In Section \ref{sec:homol-mirr-symm} we examine
Abouzaid--Smith's mirror functor (or rather, a version of it that's
adapted to our setting); we provide a description of its action on
objects that's explicit enough to enable the computations in $K_0(D
\FukExt(T^2))$ required for the proofs of the main theorems in Section
\ref{sec:isomorphism}.  Appendix \ref{sec:a_infty-triang-categ}
collects some facts on triangulated categories, cone decompositions
and Grothendieck groups; Appendix \ref{app:SESs-LSs} explains how to
get exact triangles from short exact sequences of local systems; and
finally, Appendix \ref{app:vect-bundl-ellipt} assembles a couple of
statements from Atiyah's classification of vector bundles on elliptic
curves \cite{Atiyah--VecBunEllCurve--1957}, which are used in Section
\ref{sec:homol-mirr-symm}.

\subsection*{Acknowledgements.}
\label{sec:acknowledgements}
I would like to thank my advisor, Paul Biran, for sharing so much of
his time and insight with me. I am also grateful for Octav Cornea's
interest and useful discussions, and for the anonymous referee's
helpful suggestions.

\section{The Lagrangian cobordism group}
\label{sec:lagr-cobord-group}
\subsection{Lagrangian cobordisms.}
\label{sec:Lagr-cobordisms}
We start by recalling some definitions from
\cite{Biran-Cornea--Lag-Cob-I}. For a symplectic manifold
$(M,\omega)$, we denote by $(\widetilde M, \widetilde \omega)$ the
symplectic manifold obtained by equipping $\widetilde M = \mathbb R^2
\times M$ with the split symplectic form $\widetilde \omega =
\omega_{\mathrm{std}} \oplus \omega$, where $\omega_{\mathrm{std}} =
dx \wedge dy$ is the standard symplectic form on $\mathbb R^2$. We
denote by $\pi: \widetilde M \to \mathbb R^2$ the projection to the
first factor, and given any subset $S \subset \mathbb R^2$, we write
$V\vert_S = V \cap \pi^{-1}(S)$.

We say that two ordered collections $(L_i)_{i=1}^{r}$ and
$(L_j')_{j=1}^{s}$ of Lagrangian submanifolds of $M$ are
\textit{Lagrangian cobordant} if there exists a compact cobordism
$\big(V;\coprod_i L_i, \coprod_j L_j'\big)$ together with a Lagrangian
embedding $V \to [0,1] \times \mathbb R \times M \subset \mathbb R^2
\times M$ with cylindrical ends, in the sense that there is some
$\varepsilon > 0$ such that
    \begin{equation*}
        V\vert_{[0,\varepsilon) \times \mathbb R} = \coprod_{i=1}^{r}
        ~ [0,\varepsilon) \times \{i\} \times L_i \quad \text{and} \quad V\vert_{(1-\varepsilon,1]\times \mathbb R} = \coprod_{j=1}^{s} ~
        (1-\varepsilon,1]\times \{j\} \times L_j'.
    \end{equation*}
The Lagrangian submanifold $V$ of $\widetilde M$ is called a
\textit{Lagrangian cobordism} with positive ends $(L_j')_{j=1}^s$ and
negative ends $(L_i)_{i=1}^r$. The terminology is that $V$ goes from
$(L_j')_{j=1}^s$ to $(L_i)_{i=1}^r$ and we denote this relationship by
\begin{equation*}
    V: (L_1',\dots,L_s') \leadsto (L_1,\dots,L_r).
\end{equation*}

\begin{example}
    \emph{Hamiltonian isotopy}: Let $\Phi: M \to M$ be a
    Hamiltonian diffeomorphism, and let $\phi: [0,1] \times M \to M$
    be a Hamiltonian isotopy with $\phi(0,\cdot) = \mathrm{id}$,
    $\phi(1,\cdot) = \Phi$ and which is generated by a Hamiltonian $H:
    [0,1]\times M \to \mathbb R$ such that $H(t,\cdot) \equiv 0$ for
    $t$ close to $0$ and $1$ (this condition can be achieved by suitably
    reparametrising any given Hamiltonian isotopy). Then for any
    Lagrangian submanifold $L \subset M$, the map
    \begin{equation*}
        [0,1] \times L \to \mathbb R^2 \times M, \quad (t,x) \mapsto \left(t,\, - H(t,\phi_t(x)), \,\phi_t(x)\right),
    \end{equation*}
    defines a Lagrangian cobordism $V: \phi_1(L) \leadsto L$.
\end{example}

\begin{example}
    \emph{Lagrangian surgery}: Let $L_0, L_1 \subset M$ be two
    transversely intersecting Lagrangian submanifolds. One can resolve
    the intersection points by cutting out small neighbourhoods and
    gluing in Lagrangian handles diffeomorphic to $[-1,1] \times
    S^{n-1}$ (where $n = \frac{1}{2} \mathrm{dim} \, M$). This
    produces a new Lagrangian submanifold which we denote by $L_0 \#
    L_1$. Biran--Cornea \cite{Biran-Cornea--Lag-Cob-I} show that there
    exists a Lagrangian cobordism
    \begin{equation*}
        V: L_0  \#L_1 \leadsto (L_0,L_1).
    \end{equation*}
    This construction will be of prime importance later on; we will
    describe it in more detail in Section \ref{sec:lagr-cobord-T2}.
\end{example}

\subsection{Extra data.}
Whenever Lagrangians come equipped with extra data, such as
orientations, $Spin$ or $Pin$ structures, gradings, or local systems,
it makes sense to consider cobordisms over which these data extend. We
say that two collections $(L_i^\deco)_{i=1}^r$ and
$(L_j'^\deco)_{j=1}^s$ of Lagrangians decorated with such extra data
are Lagrangian cobordant if there exists a Lagrangian cobordism $V:
(L_j')_{j=1}^s \leadsto (L_i)_{i=1}^r$ between the underlying
Lagrangians together with choices of the same types of extra data for
$V$ which restrict to the given data on the ends (provided there
exists a suitable notion of restricting to boundary components).

\subsection{The Lagrangian cobordism group.}
\label{sec:lagr-cobord-group-1}
Let $\Lags$ be the set of (suitably qualified) Lagrangian submanifolds
of $M$, and let $\langle \Lags \rangle$ be the free Abelian group
generated by $\Lags$. Denote by $R \subset \langle \Lags \rangle$ the
subgroup generated by all expressions
\begin{equation*}
    L_1 + \dots +
    L_{r} - L_1' - \dots -L_{s}' \in \langle \Lags \rangle
\end{equation*}
such that there is a (suitably qualified) Lagrangian cobordism $ V:
(L_1',\dots,L_s') \leadsto (L_1,\dots,L_r)$. The Lagrangian cobordism
group corresponding to $\Lags$ and $R$ is then defined as
\begin{equation*}
    \cobGr (M) = \langle \Lags \rangle / R,
\end{equation*}
where we suppress the dependence of $\Lags$ and $R$ in the notation. An
analogous definition applies when Lagrangians carry extra data, in
which case one imposes relations coming from cobordisms with the same
types of extra data.

As indicated, it usually makes sense to constrain the Lagrangians and
cobordisms one admits in the definition of $\cobGr(M)$, because
Lagrangian cobordism without any additional condition is a quite
flexible notion. One possibility is to require all Lagrangians and
cobordisms to be (uniformly) monotone, as is done in
\cite{Biran-Cornea--Lag-Cob-I, Biran-Cornea--Lag-Cob-II}. In the case
$M = T^2$ this paper deals with, we
will require that the Lagrangians are non-contractible curves and that
the cobordisms have vanishing Maslov class.

\begin{rmk}
    Constraining the cobordisms one allows in the definition of
    $\cobGr(M)$ obviously has an effect on what equality in
    $\cobGr(M)$ means. For example, the identity $[L] = [L']$ in the
    monotone version of $\cobGr(M)$ does not necessarily mean that
    there exists a monotone cobordism $L' \leadsto L$.  The identity
    might instead come from a monotone cobordism $(L',K) \leadsto
    (L,K)$; from this one \emph{can} create a cobordism $L' \leadsto
    L$ by connecting the ends corresponding to $K$, but in general at
    the cost of losing monotonicity.
\end{rmk}

\subsection{Additional relations from local systems.}
\label{sec:relations-from-local-systems}
Next we define a variant of the Lagrangian cobordism group that one
can define whenever the Lagrangians one studies carry local
systems. Let $\LagsExt$ be the set of all pairs $(L,E)$ where $L$ is a
Lagrangian submanifold in $M$ and where $E$ is a local system on $L$
(both with suitable qualifications depending on the context), and let
$\langle \LagsExt \rangle$ be the free Abelian group generated by
$\LagsExt$.

To define the subgroup $R^\sesLS \subset \langle \LagsExt \rangle$ of
relations we impose, consider first all expressions of the form
\begin{equation}
    \label{eq:rels-LagsExt-1}
    (L_1,E_1) + \dots + (L_r,E_r) - (L_1',E_1') - \dots - (L_s',E_s')
    \in \langle \LagsExt \rangle
\end{equation}
such that there is a Lagrangian cobordism $V: (L_1',\dots,L_s')
\leadsto (L_1,\dots,L_r)$ together with a local system $E$ on $V$
which restricts to the $E_i'$ or $E_j$ on the respective
ends. Second, consider all expressions of the form
\begin{equation}
    \label{eq:rels-LagsExt-2}
    (L,E) - (L,E') - (L,E'') \in \langle \LagsExt
    \rangle
\end{equation}
such that there exists a short exact sequence $0 \to E' \to E \to E''
\to 0$ of local systems on $L$. Now define $R^\sesLS$ to be the
subgroup of $\langle \LagsExt \rangle$ generated by all expressions of
types \eqref{eq:rels-LagsExt-1} and \eqref{eq:rels-LagsExt-2} and set
\begin{equation*}
    \cobGrExt(M) = \langle \LagsExt \rangle / R^\sesLS.
\end{equation*}
Similar to before, we suppress $\LagsExt$ and $R^\sesLS$ from the
notation, keeping in mind that the group really depends on them.

\section{The Fukaya category of $T^2$}
\label{sec:Fuk-T2}
In this section we describe the constructions of the Fukaya categories
$\Fuk(T^2)$ and $\FukExt(T^2)$ appearing in Theorems \ref{thm:MainThm}
and \ref{thm:MainThmExtended}, following essentially the general
recipe in \cite{Seidel--Fukaya-Picard-Lefschetz-2008}. The ground
field over which these $A_\infty$-categories are defined is the
Novikov field
\begin{equation*}
    \Lambda = \left\{ \sum_{i =0}^\infty c_i
        q^{a_i}~\Big|~ c_i \in \mathbb C, ~a_i \in \mathbb R, ~a_i
        < a_{i+1}, ~\lim_{i \to \infty} a_i = \infty 
    \right \}.
\end{equation*}

We will first describe the objects of $\FukExt(T^2)$ and $\Fuk(T^2)$,
and then sketch the construction of morphism spaces,
$A_\infty$-composi\-tions and derived categories.

\subsection{Objects.}
\label{sec:fukaya-category}
The objects of $\FukExt(T^2)$ are tuples $(L,\alpha,P,E)$ as follows:
\begin{itemize}
\item $L \subset T^2$ is a non-contractible simple closed curve,
\item $\alpha: L \to \mathbb R$ is a grading of $L$,
\item $P$ is a $Pin$ structure on $L$,
\item $E$ is a local system of $\Lambda$-vector spaces on $L$.
\end{itemize}
These structures are subject to certain conditions which we will
describe below. 
The pair $(\alpha,P)$ is called a \emph{brane structure} on $L$, and
the triple $(L,\alpha,P)$ is called a \emph{Lagrangian brane}. We
usually suppress the brane structure and/or the local system from the
notation.

The objects of $\Fuk(T^2)$ are Lagrangian branes $L^\deco =
(L,\alpha,P)$ with $L$, $\alpha$ and $P$ as above. We regard them as
objects of $\FukExt(T^2)$ by equipping them with trivial rank one
local systems, such as to make $\Fuk(T^2)$ a full
$A_\infty$-subcategory of $\FukExt(T^2)$. In view of that, we won't
define the morphism spaces of $\Fuk(T^2)$ separately.

\subsubsection{Gradings.}
\label{sec:gradings}
The Lagrangian Grassmannian $\mathrm{Gr}(\mathbb R^{2n})$ of $(\mathbb
R^{2n},\omega_{\mathrm{std}})$ can be naturally identified with $U(n)
/ O(n)$, and this identification induces a map $\mathrm{det}^2:
\mathrm{Gr}(\mathbb R^{2n}) \to S^1$. A \emph{grading} of a Lagrangian
subspace $\Lambda \in \mathrm{Gr}(\mathbb R^{2n})$ is a number $\alpha
\in \mathbb R$ such that $e^{2 \pi i \alpha} =
\mathrm{det}^2(\Lambda)$. (See
\cite{Seidel--Fukaya-Picard-Lefschetz-2008} for the general
definition.)

Since the tangent bundle of $T^2 = \mathbb R^2/\mathbb Z^2$ is $T^2
\times \mathbb R^2$, the Gau\ss~ map associated to a Lagrangian $L
\subset T^2$ can be viewed as a map $\Gamma_L: L \to
\mathrm{Gr}(\mathbb R^2)$. A \emph{grading} of $L$ is a continuous
function $\alpha: L \to \mathbb R$ such that $\alpha(x)$ is a grading
of $T_xL \in \mathrm{Gr}(\mathbb R^2)$ in the previous sense; that is,
$\alpha$ is a lift of the map $\mathrm{det}^2 \circ \Gamma_L: L \to S^1$
with respect to the covering $\mathbb R \to S^1$, $\alpha \mapsto e^{2
    \pi i \alpha}$.
Note that every non-contractible $L \subset T^2$ possesses a grading
because its $\Gamma_L$ is null-homotopic. Objects of $\FukExt(T^2)$
are allowed to carry all possible gradings.

A grading $\alpha$ on a curve $L \subset T^2$ induces an orientation
of $L$ because we can view $e^{\pi i \alpha} \in S^1$ as a point in
$\mathrm{Gr}^{\mathrm{or}}(\mathbb R^2)$, the Grassmannian of oriented
lines in $\mathbb R^2$; changing the grading by $\pm 1$ reverses the
induced orientation. In the following it will sometimes be convenient
to have an orientation present, and we will always equip graded curves
with this induced orientation. Conversely, given an oriented curve
$L$, we can assign a \emph{standard grading} to it as follows. Suppose
that $L$ has slope $(p,q) \in H_1(T^2;\mathbb Z)$. Then there is a
unique number $\alpha_0 \in [-1,1)$ such that
\begin{equation*}
    e^{\pi i \alpha_0}  = \frac{ p + iq}{\sqrt{p^2 + q^2}}.
\end{equation*} 
Viewed as a constant function, $\alpha_0$ is a grading for any of the
linear Lagrangians to which $L$ is isotopic. We define the standard
grading of $L$ to be the function $\alpha: L \to \mathbb R$ induced by
$\alpha_0$ via any isotopy connecting $L$ to a linear Lagrangian.

\begin{figure}[t]
    \centering
    \includegraphics[scale=0.8]{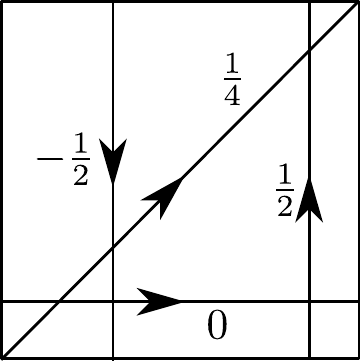}
    \caption{Oriented curves on $T^2$ and their standard grading.}
    \label{fig:std-grading}
\end{figure}

\subsubsection{Pin structures.}
\label{sec:spin-structures}
A $Pin$ structure on an rank $n$ vector bundle $F$ is a principal
$Pin_n$-bundle $P$ together with a choice of two-sheeted covering $P
\to P_{O}(F)$ which is equivariant with respect to the action of
$Pin_n$ on both sides. Here $P_{O}(F)$ is the bundle of orthonormal
frames of $F$ with respect to some auxiliary metric, which is acted
upon by $Pin_n$ via the homomorphism $Pin_n \to O_n$. A $Pin$
structure on a manifold is a $Pin$ structure on its tangent
bundle. $Pin$ structures generalize the more familiar notion of $Spin$
structures to non-oriented (even non-orientable) bundles. If bundle
$F$ is orientable, a $Pin$ structure is the same as an equivalence
class of $Spin$ structures, where two $Spin$ structures are equivalent
if they are obtained from one another by ``reversing the
orientation''. (See
\cite[Section(11i)]{Seidel--Fukaya-Picard-Lefschetz-2008} for some and
\cite{Lawson-Michelsohn--Spin-geometry} for more background on $Pin_n$
and $Pin$ structures.)

The orthonormal frame bundle of a circle $L$ is simply $L \sqcup L$,
and we can therefore think of a $Pin$ structure on $L$ as a double
cover of $L$. Hence $L$ admits precisely two different $Pin$
structures, corresponding to the trivial resp. the non-trivial
2-sheeted cover. The latter is referred to as the \textrm{bounding}
$Pin$ structure, because it is obtained by restricting the unique
$Pin$ structure on the disc to the boundary circle. The Lagrangians
appearing in the definition of $\FukExt(T^2)$ will be equipped with
the bounding $Pin$ structure. (Allowing also the other $Pin$ structure
would be possible, but would lead to some redundancy.)

\subsubsection{Local systems.}
\label{sec:local-systems}
Let $X$ be a topological space. A local system $E$ on $X$ assigns a
vector space $E_x$ to every point $x \in X$ and an isomorphism
$\pi_\gamma: E_{\gamma(0)} \to E_{\gamma(1)}$ to every path $\gamma:
[0,1] \to X$ which depends only on the homotopy class of the path
relative to the end points, in a way which is compatible with
concatenation of paths. We call $E_x$ the \emph{fibre} of $E$ at $x$,
and $\pi_\gamma$ the \emph{parallel transport} in $E$ along $\gamma$.
A local system $E$ yields a representation $\eta_{x_0}:\pi_1(X,x_0)
\to \mathrm{GL}(E_{x_0})$ for every $x_0 \in X$.  Assuming that $X$ is
path-connected, $E$ is determined up to isomorphism by one single such
representation $\eta_{x_0}$; 
in fact, $E$ is isomorphic to a local system whose fibres are all
equal to $E_{x_0}$ and with parallel transport maps constructed from
$\eta_{x_0}$. 

If the base space is an oriented circle $L$, a local system $E$ on $L$
is hence determined by specifying a single vector space $F_E$ and a
\emph{monodromy} isomorphism $M_E \in \mathrm{GL}(F)$ which encodes
the parallel transport along the preferred generator of
$\pi_1(L)$. The local systems $E$ we allow in the definition of
$\FukExt(T^2)$ are as follows: They have fibre $F_E = \Lambda^{\oplus
    n}$ for some $n \geq 1$, and all eigenvalues of the monodromy $M_E
\in \mathrm{GL}_n(\Lambda)$ have norm 1 with respect to the
non-Archimedean norm on $\Lambda$ defined by
\begin{equation*}
    \left\vert\sum_{i=0}^\infty c_i q^{a_i} \right\vert = e^{-a_0}.
\end{equation*}

\subsection{Morphisms.} Let $L^\deco_i \equiv (L_i,\alpha_i,P_i,E_i)$,
$i = 0,1$, be objects of $\FukExt(T^2)$ such that $L_0$ and $L_1$
intersect transversely. The corresponding space of \emph{Floer
    cochains} is the graded $\Lambda$-vector space with $i$-th graded
component
\begin{equation*}
    CF^i(L_0^\deco,L_1^\deco) = \bigoplus_{\substack{y \in L_0 \cap L_1\\i(y) \,=\,            i}} \mathrm{Hom}(E_{0,y}, E_{1,y}),
\end{equation*}
where $\mathrm{Hom}(E_{0,y}, E_{1,y})$ is the space of homomorphisms
between the fibres of the local systems at $y \in L_0 \cap L_1$, and
$i(y)$ is the \emph{index} of $y$, defined by
\begin{equation*}
    i(y) \equiv i(y;L_0,L_1) = \lfloor \alpha_1(y) - \alpha_0(y)
    \rfloor + 1,
\end{equation*}
where $\lfloor \cdot \rfloor$ is the next lowest integer.

Let $Z = \mathbb R \times [0,1]$ be the strip with coordinates $(s,t)$
and equipped with the usual complex structure $j_Z$. Given $y_0,y_1
\in L_0 \cap L_1$, denote by $\mathcal M_Z(y_0,y_1)$ the space of all
maps $u: Z \to T^2$ satisfying
\begin{equation}
    \label{eqn:Floer-eqn}
    \partial_su + J(t,u) \partial_t u = 0
\end{equation}
for a generic $\omega$-compatible $t$-dependent almost complex
structure $J$ on $T^2$, and with boundary and asymptotic conditions
given by
\begin{equation*}
    u(s,0) \in L_0, ~ u(s,1) \in L_1, ~ \lim_{s \to
        -\infty}u(s,\cdot) = y_0, ~ \lim_{s \to +\infty} u(s,\cdot) = y_1.
\end{equation*}
$\mathcal M_Z(y_0,y_1)$ carries a natural $\mathbb R$-action given by
translation in the $\mathbb R$-variable. We denote by $ \mathcal
M^{1+1}(y_0,y_1) = \mathcal M_Z(y_0,y_1) / \mathbb R$ the moduli space
obtained by quotienting out this action.

The Floer differential $\partial: CF(L_0^\deco,L_1^\deco) \to
CF(L_0^\deco,L_1^\deco)[1]$ is then defined on generators $\phi_1 \in
\mathrm{Hom}(E_{0,y_1},E_{1,y_1})$ of $CF(L_0^\deco,L_1^\deco)$ by 
\begin{equation*}
    \partial \phi_1 = (-1)^{i(y_1)}\bigoplus_{y_0 \in L_0 \cap L_1} \sum_{u} \mathrm{sgn}(u) ~\pi_1^u \circ \phi \circ \pi_0^u ~q^{\omega(u)},
\end{equation*}
where the second sum runs over the zero-dimensional component $\mathcal
M^{1+1}(y_0,y_1)^0$ of the moduli space (which is a discrete
set). Here $\pi_i^u$, for $i = 0,1$, denotes parallel transport in the
local system $E_i$ along the boundary component $u(\mathbb R \times
\{i\}) \subset L_i$ of the strip $u(Z)$, and $q^{\omega(u)}$ is an
element of $\Lambda$ that encodes the symplectic area of the strip
$u$. Finally, $\mathrm{sgn}(u) \in \{\pm 1\}$ is a sign whose
determination we will describe in Section \ref{sec:signs}.

The definition of $CF(L_0^\deco,L_1^\deco)$ for non-transversely
intersecting $L_0,$ $L_1$ requires the use of Hamiltonian
perturbations: One fixes a Floer datum $(H,J)$ for every such pair,
consisting of a Hamiltonian function $H$ such that $\phi_1^H(L_0)
\pitchfork L_1$ and an almost complex structure $J$, and then
considers an analogue of equation \eqref{eqn:Floer-eqn} with an
additional perturbation term on the right-hand side.

In either case, the space of morphisms in $\FukExt(T^2)$ from
$L_0^\deco$ to $L_1^\deco$ is defined as a graded vector space over
$\Lambda$ by
\begin{equation*}
    \mathrm{hom}(L_0^\deco,L_1^\deco) = CF(L_0^\deco,L_1^\deco),
\end{equation*}
and the $A_\infty$-structure map $\mu^1: \mathrm{hom}(L_0,L_1) \to
\mathrm{hom}(L_0,L_1)[1]$ of order one is the Floer differential
$\partial$.

\subsection{$A_\infty$-compositions.}
\label{sec:a_infty-compositions}
We give a brief description of the higher
$A_\infty$-com\-po\-si\-tions
\begin{equation*}
    \mu^d: \mathrm{hom}(L_{d-1}^\deco,L_d^\deco) \otimes \cdots \otimes
    \mathrm{hom}(L_{0}^\deco,L_1^\deco) \to
    \mathrm{hom}(L_{0}^\deco,L_d^\deco)[2-d], ~ d \geq 2,
\end{equation*}
again limiting ourselves to the case of mutually transverse $L_i$. We
refer to \cite{Seidel--Fukaya-Picard-Lefschetz-2008} for the proof
that these really define an $A_\infty$-structure.

\begin{figure}[t]
    \centering
    \includegraphics[scale=1]{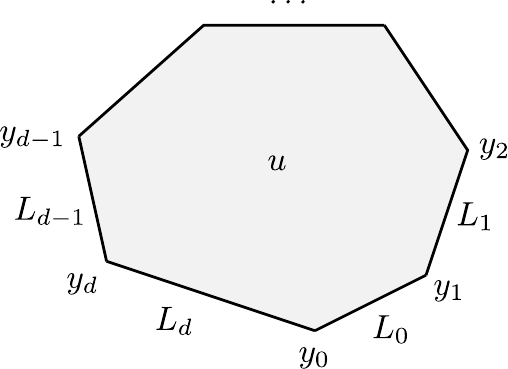}
    \caption{A pseudo-holomorphic polygon.}
    \label{fig:pseudohol-polygon}
\end{figure}

Let $L_0^\deco,\dots,L_d^\deco$ be objects of $\FukExt(T^2)$, and let
$y_0 \in L_0 \cap L_d$ and $y_i \in L_{i-1} \cap L_i$, for $i =
1,\dots,d$. Moreover, let $S$ be a disc with one incoming and
$d$ outgoing boundary punctures. The moduli space of $\mathcal
M^{d+1}(y_0,\dots,y_d)$ of \emph{pseudo-holomorphic polygons}
associated to this collection Lagrangians and points is the set of all
maps $u \in C^\infty(S,T^2)$ solving the equation
\begin{equation*}
    Du(z)  + J(z,u) \circ Du(z) \circ j_S = 0,
\end{equation*}
where $J$ is a generic $\omega$-compatible almost complex structure
depending on $z \in S$, with boundary conditions given by the $L_i$
and asymptotic conditions at the punctures given by the $y_i$ as
indicated in Figure \ref{fig:pseudohol-polygon}. The precise general
definitions require again choices of Floer data and additional choices
of perturbation data which lead to the appearance of an inhomogeneity
on the righthand side of the equation; see \cite[Section
8(f)]{Seidel--Fukaya-Picard-Lefschetz-2008} for details.

Given now $\phi_1,\dots,\phi_d$ with $\phi_i \in
\mathrm{Hom}(E_{i-1,y_i},E_{i,y_i})$, the corresponding output of
$\mu^d$ is
\begin{equation*}
    \mu^d(\phi_d,\dots,\phi_1) = (-1)^\dagger\bigoplus_{y_0 \in L_0 \cap L_d}
    \sum_{u} \mathrm{sgn}(u)  ~\pi_d^u \circ
    \phi_d \circ \pi_{d-1}^u \circ \dots \circ \phi_1 \circ \pi_0^u~q^{\omega(u)},
\end{equation*}
which is an element of $\mathrm{Hom}(E_{0,y_0},E_{d,y_0}) \subset
CF(L_0^\deco,L_d^\deco)$. Here the second sum runs over all elements
$u \in \mathcal M^{d+1}(y_0,\dots,y_d)^0$, the zero-dimensional
component of the moduli space. $\pi^u_i$ is the parallel transport in
$E_i$ along the boundary component of $u(S)$ that lies in $L_i$,
and $q^{\omega(u)} \in \Lambda$ encodes the symplectic area of
$u$. The determination of the sign $\mathrm{sgn}(u) \in \{ \pm 1\}$
will be explained below; $(-1)^\dagger$ is an additional sign, with
$\dagger = \sum_{k=1}^d k i(y_k)$. 

\begin{rmk}
    The set $\mathcal M^{d+1}(y_0,\dots,y_d)^0$ is \emph{infinite} in
    general (cf. \cite{polishchuk-zaslow--HMS-elliptic-curve99} for
    the case $d = 2$). Nevertheless, the sums above converge in
    $\mathrm{Hom}(E_{0,y_0},E_{d,y_0})$ with respect to the topology
    induced by the non-Ar\-chi\-me\-de\-an norm on $\Lambda$. This
    follows from two facts: First, the monodromies of the $E_i$ 
    have norm 1; second, the areas $\omega(u)$ of the polygons
    entering the count tend to $\infty$, as one can see by thinking of
    their lifts to the universal cover of $T^2$.
\end{rmk}
\subsection{Signs.}
\label{sec:signs} 
Let $L_0^\deco, L_1^\deco$ be Lagrangian branes intersecting
transversely at some point $y \in L_0 \cap L_1$. One can associate to
$y$, considered as a morphism from $L_0^\deco$ to $L_1^\deco$, a
1-dimensional real vector space $o(y)$ called the \emph{orientation
    space} of $y$. We refer to
\cite{Seidel--Fukaya-Picard-Lefschetz-2008} for the precise definition
(also cf. Section \ref{sec:curves-in-M-vs-Mtilde} below). At this
point we only remark that $o(y)$ depends on a choice of homotopy class
of paths from $T_yL_0$ to $T_yL_1$ in the Lagrangian Grassmannian;
this homotopy class in turn is canonically determined by the brane
structures.

Given Lagrangian branes $L_0^\deco,\dots, L_d^\deco$ and intersection
points $y_i$ as in Section \ref{sec:a_infty-compositions}, there
exists a preferred isomorphism
\begin{equation*}
    \label{eq:or-of-mod-spaces}
    \Lambda^{\mathrm{top}}(T_u\mathcal M^{d+1}(y_0,\dots,y_d)) \cong o(y_0)
    \otimes o(y_1)^\vee \otimes \dots \otimes o(y_d)^\vee
\end{equation*}
for every regular $u \in \mathcal M^{d+1}(y_0,\dots,y_d)$, see
\cite[Section (12b)]{Seidel--Fukaya-Picard-Lefschetz-2008}.  In
particular, whenever $u$ is isolated (that is, an element of the
zero-dimensional component of the moduli space), this isomorphism
yields a preferred element
\begin{equation*}
    c_u \in o(y_0)\otimes o(y_1)^\vee \otimes
    \dots \otimes o(y_d)^\vee,
\end{equation*}
because $\Lambda^{\mathrm{top}}(T_u\mathcal M^{d+1}(y_0,\dots,y_d)) =
\mathbb R$ for isolated $u$.

The signs $\mathrm{sgn}(u)$ are then defined as follows: We choose,
arbitrarily and once and for all, an orientation $\mathfrak o_{y_i}$
of $o(y_i)$ for every $y_i$, which induces an orientation for every
$o(y_0) \otimes o(y_1)^\vee \otimes \dots \otimes o(y_d)^\vee$. Then
we set $\mathrm{sgn}(u) = \pm 1$ according to whether $c_u$ is
positive or negative with respect to this orientation.

\begin{rmk}
    \label{rmk:signs-here-vs-Seidel}
    An equivalent way of dealing with signs would be to adopt Seidel's
    basis-free approach from
    \cite{Seidel--Fukaya-Picard-Lefschetz-2008} and to define Floer
    complexes as $CF(L_0^\deco,L_1^\deco) = \bigoplus_y \vert o(y)
    \vert_\Lambda \otimes \mathrm{Hom}(E_{0,y},E_{1,y})$, where $\vert
    o(y)\vert_\Lambda$ is the \emph{$\Lambda$-nor\-mal\-i\-za\-tion}
    of $o(y)$; that is, $\vert o(y)\vert_\Lambda$ is the 1-dimensional
    $\Lambda$-vector space obtained by taking the vector space
    generated by the two orientations of $o(y)$, and imposing the
    relation that the two generators add up to zero. The
    $A_\infty$-compositions would then defined by combining the
    preferred elements $\vert c_u \vert_\Lambda \in \vert
    o(y_0)\vert_\Lambda \otimes \vert o(y_1)\vert_\Lambda^\vee \otimes
    \dots \otimes \vert o(y_d)\vert_\Lambda^\vee$ induced by the $c_u$
    with the parallel transport maps $\pi^u_i$. The approach to signs
    we have chosen is slightly more geometric (though a bit less
    elegant), since the generators of Floer complexes are actual
    homomorphisms between fibres of local systems. To translate
    between the two approaches, one identifies the two versions of the
    Floer complexes by mapping a generator $\phi \in
    \mathrm{Hom}(E_{0,y},E_{1,y})$ of the first version to $[\mathfrak
    o_y] \otimes \phi \in \vert o(y) \vert_{\Lambda} \otimes
    \mathrm{Hom}(E_{0,y},E_{1,y})$, which is a generator of the second
    version; here $\mathfrak o_y$ is the chosen orientation for
    $o(y)$. This intertwines the respective versions of the
    $A_\infty$-compositions.
\end{rmk}

\subsection{The derived Fukaya category.} The derived Fukaya category
$D \FukExt(T^2)$ is constructed from $\FukExt(T^2)$ by a purely
algebraic procedure. One first completes it to a triangulated
$A_\infty$-category $\overline{\FukExt}(T^2)$, and then sets $D
\FukExt(T^2) = H^0(\overline {\FukExt}(T^2))$. A further completion
that formally introduces images for all idempotent morphisms yields
the split-closed derived Fukaya category $D^\pi \FukExt(T^2)$.

These categories are triangulated categories in the classical sense,
meaning that they possess a shift functor and a class of exact
triangles (see Appendix \ref{sec:a_infty-triang-categ} for a brief
description of what that means). The shift functor is realized
geometrically by changing the brane structure in a certain way. In our
case its effect is simply a shift of the grading by 1,
\begin{equation*}
    (L,\alpha,P,E)[1] = (L,\alpha-1,P,E).
\end{equation*}
(See \cite[(11k)]{Seidel--Fukaya-Picard-Lefschetz-2008}; the $Pin$
structure remains unchanged since we are dealing with curves, which
have trivial tangent bundle.)

\section{Cobordisms and cone decompositions}
\label{sec:lagr-cobord-cones}
Let $(M,\omega)$ be a symplectic manifold and suppose that
$\FukExt(M)$ is defined as a graded $A_\infty$-category over the
Novikov field $\Lambda$, with objects Lagrangian branes carrying local
systems.

\begin{thm}
    \label{thm:cob-cone-decomp}
    Let $L^\deco, L_1^\deco,\dots,L_r^\deco \in \FukExt(M)$ and
    suppose that there exists a Lagrangian cobordism $V^\deco: L^\deco
    \leadsto (L_1^\deco,\dots,L_k^\deco)$.  Then $L^\deco$ admits an
    iterated cone decomposition in $D \FukExt(M)$ with
    linearization $(L_1^\deco,\dots,L_k^\deco)$.
\end{thm}

This holds provided that certain technical conditions are satisfied,
which will be addressed in Section \ref{sec:differences-setting}. We
suppress brane structures and local systems from the notation, i.e.,
we write $L$ instead of $(L,\alpha,P,E)$. By a Lagrangian cobordism
between such objects we mean a Lagrangian cobordism between the
underlying Lagrangians which is equipped with the same type of extra
data in a compatible way. The definitions of \emph{iterated cone
    decomposition} and \emph{linearization} are given in Appendix
\ref{sec:a_infty-triang-categ}.

Theorem \ref{thm:cob-cone-decomp}, which will follow from Proposition
\ref{prop:cob-cone-decomp} below, which is an adaptation of an
analogous statement in \cite{Biran-Cornea--Lag-Cob-II} to the present
setting, and this section serves to explain the (small) modifications
of Biran-Cornea's arguments required for its proof.

\subsection{Differences to the setting in
    \cite{Biran-Cornea--Lag-Cob-II}.}
\label{sec:differences-setting}
The version of the Fukaya category considered in
\cite{Biran-Cornea--Lag-Cob-II} is linear over $\mathbb Z/2\mathbb Z$
and ungraded, and has as objects plain Lagrangians $L$, without brane
structures or local systems. Our main point will be to explain the
inclusion of gradings and signs in the proofs given in
\cite{Biran-Cornea--Lag-Cob-II}; including the local systems is
standard and doesn't pose any additional difficulties. Some small
differences are due to the fact that we work with Floer cohomology (in
contrast to Floer homology).

Biran-Cornea require that $(M,\omega)$, as well as all Lagrangians $L$
and cobordisms $V$ involved, are uniformly monotone, and that the maps
\begin{equation*}
    \pi_1(L) \to \pi_1(M) \quad \text{and} \quad \pi_1(V) \to \pi_1(\widetilde M)
\end{equation*}
induced by the inclusions vanish for all of them. These conditions are
needed to prove compactness of moduli spaces of pseudo-holomorphic
curves.  More precisely, the conditions on fundamental groups are used
to obtain bounds on the areas of such curves; 
we do not have to impose these conditions, because we use the Novikov
field $\Lambda$ to encode series of curves with areas tending to
infinity. The monotonicity assumption is needed to rule out the
bubbling off of pseudo-holomorphic discs.

Discs bubbles can, however, also be excluded in some situations which
do not fit into the monotone setting. For example, we will need
Theorem \ref{thm:cob-cone-decomp} in the case that $L, L_1,\dots,L_r$
are non-contractible curves in $T^2$ and $V$ is a Lagrangian cobordism
coming from iterated surgery of such curves, as described in Section
\ref{sec:lagr-cobord-T2}. Our curves clearly can't bound any discs as
they are non-contractible. Moreover, our cobordisms will be shown to
have vanishing Maslov class (Proposition \ref{lem:Gauss-cobsT2}); as
the expected dimension of the moduli space of Maslov zero
pseudo-holomorphic discs with boundary on a Lagrangian in a symplectic
4-manifold is $-1$, such discs exist only for almost complex
structures $J$ belonging to a codimension one stratum $\mathcal J_0
\subset \mathcal J$ of the space of compatible almost complex
structures. Thus bubbling can be excluded as long as long as one only
works with almost complex structure $J$ in some fixed component of
$\mathcal J \setminus \mathcal J_0$. (This also means that the
resulting objects are only invariant with respect to changing $J$
within this component, but that is not a problem for our specific
purposes.)

For the arguments of this section, it will not matter which precise
mechanism prevents disc bubbling. We will simply assume that we are in
one of these favourable settings.

\subsection{The categories $\mathscr B_{\gamma,h}$.}
\label{sec:auxiliary-categories}
As in \cite{Biran-Cornea--Lag-Cob-II}, we consider
$A_\infty$-cat\-e\-go\-ries $\mathscr B_{\gamma,h}$ associated to
pairs $(\gamma,h)$ as follows: First, $\gamma \subset \mathbb R^2$ is
a properly embedded curve diffeomorphic to $\mathbb R$ with
``horizontal ends'' whose $y$-coordinates are in $\frac{1}{2}
\,\mathbb Z$; second, $h: \mathbb R^2 \to \mathbb R$ is a
\emph{profile function}, whose associated \emph{extended profile
    function} we denote by $h':\mathbb R^2 \to \mathbb R$.

\begin{figure}[t]
    \centering
    \includegraphics[scale=1.4]{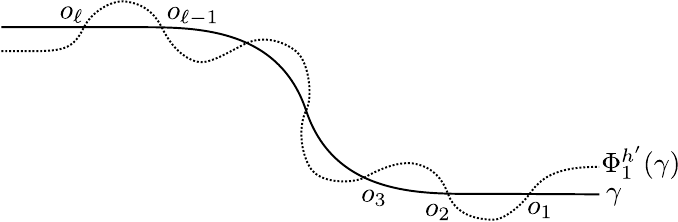}
    \caption{The curves $\gamma$ and $\Phi_1^{h'}(\gamma)$ used to
        define $\mathscr B_{\gamma,h}$.}
    \label{fig:gamma-curves}
\end{figure}

These terms are defined in Sections 3.2 and 4.1 of
\cite{Biran-Cornea--Lag-Cob-II}. What is really only important for us,
is that $h$ and $h'$ are assumed to be such that $\Phi^{h'}_1(\gamma)$
looks as in Figure \ref{fig:gamma-curves}, where $\Phi^{h'}_1$ is the
time one map of the Hamiltonian flow of $h'$, and $\gamma$ is a
typical specimen of the type of curves considered. That is, there
should be an \emph{odd} number of intersection points
$o_1,\dots,o_\ell \in \gamma \cap \Phi_1^{h'}(\gamma)$, and
$\Phi_1^{h'}(\gamma)$ should have horizontal ends such that its
negative end lies below the negative end of $\gamma$.

\begin{rmk}
    Note that in \cite{Biran-Cornea--Lag-Cob-II} the requirement for
    $h'$ is that $(\Phi_1^{h'})^{-1}(\gamma)$ looks as in Figure
    \ref{fig:gamma-curves}. The difference is due to our use of
    \emph{co}homology. 
\end{rmk}

The \emph{objects} of $\mathscr B_{\gamma,h}$ are pairs $(\widetilde L^\deco,
\widetilde E)$, where $\widetilde L^\deco$ is a Lagrangian brane in
$\widetilde M = \mathbb R^2 \times M$ whose underlying Lagrangian is
of the form $\widetilde L = \gamma \times L$ for some $L^\deco \in
\mathrm{Ob}\,\FukExt(M)$, and $\widetilde E$ is a local system on
$\widetilde L$. 
We usually abbreviate $\widetilde L^\deco \equiv (\widetilde
L^\deco,\widetilde E)$.

The space of \emph{morphisms} in $\mathscr B_{\gamma,h}$ between two
such $\widetilde L_0^\deco$ and $\widetilde L_1^\deco$ is
\begin{equation*}
    \mathrm{hom}_{\mathscr B_{\gamma,h}}(\widetilde L_0^\deco,
    \widetilde L_1^\deco) = CF(\widetilde L_0^\deco,\widetilde L_1^\deco),
\end{equation*}
where the Floer complex on the right-hand side is defined with respect
to a Floer datum $(\widetilde H, \widetilde J)$ for $(\widetilde L_0,
\widetilde L_1)$ which is chosen compatibly with the Floer datum
$(H,J)$ for $(L_0,L_1)$ used to define the Floer complex
$CF(L_0^\deco,L_1^\deco) =
\mathrm{hom}_{\FukExt(M)}(L_0^\deco,L_1^\deco)$.

Concretely, the first part of the Floer datum is the time-dependent
Hamiltonian function $\widetilde H = h' \oplus H$. With this choice,
the intersection points of $\Phi_1^{\widetilde H}(\widetilde L_0)$ and
$\widetilde L_1$ are of the form $(o_j,y)$, with $o_j \in
\Phi_1^{h'}(\gamma) \cap \gamma$ and $y \in \Phi_1^H(L_0) \cap
L_1$. Hence $CF(\widetilde L^\deco, \widetilde L'^\deco)$ is generated
by elements of the form
\begin{equation*}
    \phi^{o_j} \in \mathrm{Hom}(\widetilde
    E_{0,(o_j,y)},\widetilde E_{1,(o_j,y)}),
\end{equation*}
$j= 1,\dots,\ell$, and there is an obvious splitting of vector
spaces
\begin{equation}
    \label{eq:splitting-of-auxiliary-homs}
    CF(\widetilde L_0^\deco, \widetilde L_1^\deco) =
    \bigoplus_{j=1}^\ell CF(\widetilde L_0^\deco, \widetilde L_1^\deco)^{o_j},
\end{equation}
where $CF(L_0^\deco,L_1^\deco)^{o_j}$ is the summand consisting of
homomorphisms between the fibres of the local systems over
$(o_j,y)$. The second part of the Floer datum is the time-dependent
almost complex structure $\widetilde J(t) = i_{h'}(t) \oplus J(t)$,
with $i_{h'}(t) = (\phi_t^{h'})_*i$, where $i$ is the standard complex
structure on $\mathbb R^2$. 

With these choices of Floer data and further choices of perturbation
data, one constructs the $A_\infty$-compositions $\mu_{\mathscr
    B_{\gamma,h}}^d$ by combining the description in
\cite{Biran-Cornea--Lag-Cob-II} with obvious modifications due to the
presence of local systems and the use of cohomology, and building in
signs as described in Section \ref{sec:signs}.

\subsection{The functors $c_j: \mathscr B_{\gamma,h} \to \FukExt(M)$.}
\label{sec:quasi-isom-mathscr}
For every odd $j$ with $1\leq j \leq \ell$, there exists a natural
$A_\infty$-functor
\begin{equation*}
    c_j \equiv c_{\gamma,h,j}: \mathscr B_{\gamma,h} \to \FukExt(M),
\end{equation*}
which will turn out to be a quasi-isomorphism. We start by describing
the action of these functors on objects; their full definition will be
given in Section \ref{sec:a_infty-quasi-isom}.

Let $(\widetilde L^\deco, \widetilde E)$ be an object of $\mathscr
B_{\gamma,h}$, and denote by $\widetilde \alpha: \widetilde L \to
\mathbb R$ the grading and by $\widetilde P$ the $Pin$ structure which
together form the brane structure on $\widetilde L$.  Recall that
$\widetilde L = \gamma \times L$ for some Lagrangian $L \subset
M$. The image of $(\widetilde L^\deco,\widetilde E)$ under $c_j$, for
all $j$, is
\begin{equation*}
    c_j (\widetilde L^\deco,\widetilde E) = (L^\deco,E),
\end{equation*}
where $L^\deco$ and $E$ obtained from $\widetilde L^\deco$ and
$\widetilde E$ as follows. Consider the inclusion $L \hookrightarrow
\widetilde L$ of $L$ as a fibre of $\widetilde L \to \gamma$ over the
negative horizontal end of $\gamma$. Then the grading, the $Pin$
structure, and the local system on $L$ are the pullbacks of the
corresponding data on $\widetilde L$ via this inclusion. With the
definitions of gradings, $Pin$ structures and local systems given in
Section \ref{sec:fukaya-category}, it should be fairly obvious what
pulling back means. As for pulling back a $Pin$ structure $\widetilde
S$ from $\widetilde L$, we include $P_{O}(TL) \hookrightarrow
P_{O}(T\widetilde L)$ by completing frames in $TL$ to frames in $T
\widetilde L$ using the right-pointing tangent vector to $\gamma$;
then we restrict $\widetilde P$ accordingly, thinking of it as a
double cover of $P_{O}(T\widetilde L)$.

\subsection{Curves in $M$ vs. curves in $\widetilde M$.}
\label{sec:curves-in-M-vs-Mtilde}
To relate the $A_\infty$-composi\-tions in $\mathscr B_{\gamma,h}$ in
terms with those of $\FukExt(M)$, we need to add to the discussion in
\cite{Biran-Cornea--Lag-Cob-II} a verification that $\widetilde
J$-holomorphic curves $\widetilde u$ in $ \widetilde M$ living in a
fibre of $\widetilde M \to \mathbb R^2$ carry the same sign as the
corresponding $J$-holomorphic curves $u$ in $M$.

We need a bit of preparation for that. Take objects $\widetilde
L_0^\deco$ and $\widetilde L_1^\deco$ of $\mathscr B_{\gamma,h}$ and
let $L_0^\deco = c_j(\widetilde L_0^\deco)$ and $L_1^\deco=
c_j(\widetilde L_1^\deco)$ be the corresponding objects of
$\FukExt(M)$. Let
\begin{equation*}
    \widetilde y = (o_j,y) \in \Phi_1^{\widetilde H}(\widetilde L_0)
    \cap L_1,
\end{equation*}
where $o_j \in \Phi_1^{h'}(\gamma) \cap \gamma$ for some $j =
1,\dots,\ell$, and $y \in \Phi_1^H(L_0) \cap L_1$. Here $(\widetilde
H,\widetilde J)$ and $(H,J)$ are Floer data corresponding to
$(\widetilde L_0, \widetilde L_1)$ and $(L_0,L_1)$ which are related
as described in Section \ref{sec:auxiliary-categories}.

We now associate an index and an orientation space to the intersection
points $y$ and $\widetilde y$ as described in \cite[Section
11]{Seidel--Fukaya-Picard-Lefschetz-2008}. To explain this for $y$,
consider the Lagrangian subspaces
\begin{equation*}
    \Lambda_{0} =
    T_y \Phi^H_1( L_0)\quad  \text{and} \quad \Lambda_{1} = T_y L_1
\end{equation*}
of $T_yM$; these spaces carry natural brane structures induced by
those on $L_0$ and $L_1$. Choose a generic path
\begin{equation*}
    \lambda \in  \Omega^-(\mathrm{Gr}(T_yM);\Lambda_0,\Lambda_1),
\end{equation*}
such that the Lagrangian subbundle of $[0,1] \times T_yM$ obtained
from $\lambda$ admits a grading which restricts to the given gradings
of $\Lambda_0$ and $\Lambda_1$ on the ends; here
$\Omega^-(\mathrm{Gr}(T_yM);\Lambda_0,\Lambda_1)$ is the space of all
paths from $\Lambda_0$ to $\Lambda_1$ in $\mathrm{Gr}(T_yM)$ such that
the \emph{crossing form} $q_\lambda(s)$ is negative definite at $s =
1$ (see \cite[Section (11f)]{Seidel--Fukaya-Picard-Lefschetz-2008} for
the definition). The index and orientation space of $y$ are then
defined as
\begin{equation*}
    i(y) = \sum_{s < 1} \mathrm{sign}(q_\lambda(s))\quad \text{and}\quad
    o(y) = \bigotimes_{s < 1} (\lambda(s) \cap
    \lambda(1))^{\mathrm{sign}(q_\lambda(s))}.
\end{equation*}
The genericity assumption on $\lambda$ implies that the intersection
$\lambda(s) \cap \lambda(1)$, $s < 1$, is $\{0\}$ except at finitely
many points, where it is a one-dimensional real vector space. Thus the
definitions make sense, and $o(y)$ is a one-dimensional real vector
space. The definitions of $i(\widetilde y)$ and $o(\widetilde y)$ are
analogous.

\begin{lem}
    \label{lem:comp-of-indices-and-orientations}
    Given $y$ and $\widetilde y = (o_j,y)$ as above, the indices
    satisfy $i(\widetilde y) = i(y)$ if $j$ is odd, and $i(\widetilde
    y) = i(y) + 1$ if $j$ is even. Moreover, there are canonical
    isomorphisms $o(\widetilde y) \cong o(y)$ if $j$ is odd, and
    $o(\widetilde y) \cong T_{o_j}\gamma \otimes o(y)$ if $j$ is even.
\end{lem}
\begin{proof}
    Set $\ell_0 = T_{o_j}(\Phi_1^{h'}(\gamma))$ and $\ell_1 =
    T_{o_j}\gamma$. Choose two generic paths $\lambda \in
    \Omega^-(\mathrm{Gr}(\mathbb R^2);\ell_0, \ell_1)$ and $\lambda'
    \in \Omega^-(\mathrm{Gr}(T_yM);\Lambda_0,\Lambda_1)$ which satisfy
    the required compatibility with the gradings. 
    \begin{figure}[h]
        \centering
        \includegraphics[scale=1.12]{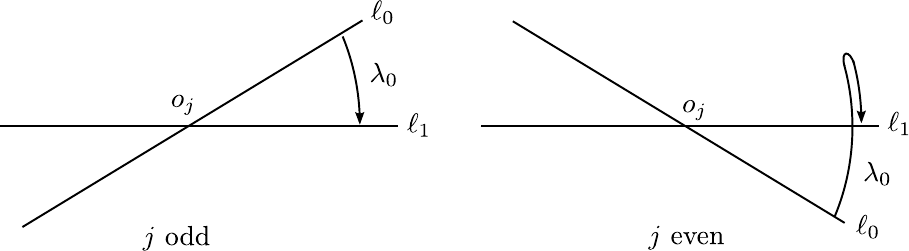}
        \caption{The paths $\lambda \in \Omega^-(\mathrm{Gr}(\mathbb
            R^2);\ell_0,\ell_1)$.}
        \label{fig:lambda}
    \end{figure}
    For $\lambda$ this means that it looks as shown in Figure
    \ref{fig:lambda}, depending on the parity of $j$ (the negative
    definiteness condition on the crossing form means $\lambda$ must
    approach $\ell_1$ from above at $s = 1$). Then the path
    $\widetilde \lambda = (\lambda, \lambda'): [0,1] \to
    \mathrm{Gr}(T_{\widetilde y}\widetilde M)$ lies in
    $\Omega^-(\mathrm{Gr}(T_{\widetilde y}\widetilde M); \widetilde
    \Lambda_0, \widetilde \Lambda_1)$ and is also compatible with the
    gradings, as a consequence of how the gradings on the $L_i$ and
    the $\widetilde L_i$ are related. The claimed statements now
    follows immediately from the definitions of the indices and of the
    orientation spaces.
\end{proof}

Recall that the moduli spaces $\mathcal M^{d+1}(y_0,\dots, y_d)$
appearing in the definition of the $A_\infty$-compositions of
$\FukExt(M)$ are spaces of maps $u:S \to M$ satisfying a
Cauchy-Riemann type equation, where $S$ is a $(d+1)$-pointed disc.
One commonly views $\mathcal M^{d+1}(y_0,\dots, y_d)$ as the zero set
of a section of $\mathcal E \to \mathcal B$, where $\mathcal B$ is a
Banach manifold of maps $u: S \to M$ which are locally of class
$W^{1,p}$ and satisfy appropriate boundary and asymptotic conditions,
and where $\mathcal E \to \mathcal B$ is the vector bundle whose fibre
at $u$ is $L^p(S;\Lambda^{0,1}T^*S \otimes E)$, with $E = u^*TM$. The
linearization of the defining section at $u \in \mathcal
M^{d+1}(y_0,\dots, y_d)$ is a non-degenerate Cauchy-Riemann operator
(in the sense of \cite[Section
(8h)]{Seidel--Fukaya-Picard-Lefschetz-2008})
\begin{equation*}
    D_u: W^{1,p}(S;E,F) \to L^p(S;\Lambda^{0,1}T^*S \otimes E),
\end{equation*}
where $F \subset E \vert_{ \partial S}$ is the Lagrangian subbundle
over $\partial S$ induced by the Lagrangians $L_i$ to which $u$ maps
the components of $\partial S$. 
Analogous statements apply to the moduli spaces $\mathcal
M^{d+1}(\widetilde y_0,\dots,\widetilde y_d)$ used to define $\mathscr
B_{\gamma,h}$. 

Suppose now that choices of Floer data have been made as in Section
\ref{sec:quasi-isom-mathscr}, and choices of perturbation data as in
\cite[Section 4.2]{Biran-Cornea--Lag-Cob-II}. Let $o_j \in
\Phi_1^{h'}(\gamma) \cap \gamma$ for some \emph{odd} $j$. Then for
every curve $u \in \mathcal M^{d+1}(y_0,\dots, y_d)$, the curve
$\widetilde u = (o_j,u)$ lies in $\mathcal M^{d+1}(\widetilde
y_0,\dots,\widetilde y_d)$ (after applying a small perturbation which
only affects the $\mathbb C$-compo\-nent of the curve -- in
\cite{Biran-Cornea--Lag-Cob-II} this is called a naturality
transformation). 
This is a consequence of the choice of perturbation datum used in the
equation defining $\mathcal M^{d+1}(\widetilde y_0,\dots,\widetilde
y_d)$, which around the odd points splits into a planar part and a
vertical part that is identical to the perturbation datum appearing in
the equation defining $\mathcal M^{d+1}(y_0,\dots, y_d)$. On the level
of linearized operators, we have a splitting
\begin{equation*}
    D_{\widetilde u} = D_{\widetilde u,\mathbb C} \oplus D_u
\end{equation*}
with respect to the canonical identifications
\begin{equation*}
    \begin{aligned}
        W^{1,p}(S;\widetilde E, \widetilde F) ~ &\cong
        ~W^{1,p}(S;\mathbb
        C, T_{o_j}\gamma) \oplus W^{1,p}(S;E,F)\\
        L^p(S;\Lambda^{0,1}T^*S \otimes \widetilde E) ~&\cong~
        L^p(\Lambda^{0,1}T^*S \otimes \mathbb C) \oplus
        L^p(\Lambda^{0,1}T^*S \otimes E)
    \end{aligned}
\end{equation*}
of the relevant spaces of sections (here $E = u^*TM$, $\widetilde E =
\widetilde u^*T \widetilde M$, and $\widetilde F$, $F$ are the
Lagrangian subbundles over $\partial S$ corresponding to the boundary
conditions for $\widetilde u$ and $u$). Again as a consequence of the
form the Floer and perturbation data around the $o_j$ with $j$ odd,
the planar operator $D_{\widetilde u, \mathbb C}$ has
$\mathrm{index}\,D_{\widetilde u,\mathbb C} = 0$,
$\mathrm{ker}\,D_{\widetilde u, \mathbb C} = \{0\}$ and
$\mathrm{coker}\,D_{\widetilde u,\mathbb C} = \{0\}$. Therefore, there
are canonical isomorphisms
\begin{equation*}
    \mathrm{ker}\,D_{\widetilde u} ~\cong~ \mathrm{ker}\,D_u,
\end{equation*}
and moreover the map $u \mapsto \widetilde u$ identifies
\begin{equation*}
    \mathcal M^{d+1} (y_0,\dots,y_d) \cong \mathcal M^{d+1}(\widetilde
    y_0,\dots,\widetilde y_d).
\end{equation*}

\begin{rmk}
    In fact, if the number of outgoing points is $d = 1$, all these
    statements are also true in case $j$ is even, because the index of
    the planar operator $D_{\widetilde u,\mathbb C}$ is zero in this
    case.  In contrast, this index is negative for even $j$ if $d$ is
    $\geq 2$ (cf. the relation between $i(y)$ and $i(\widetilde y)$ in
    Lemma \ref{lem:comp-of-indices-and-orientations} and the index
    formula for operators in
    \cite[Prop. 11.13]{Seidel--Fukaya-Picard-Lefschetz-2008}).
\end{rmk}

\subsection{Compatible choices of signs.}
\label{sec:compatible-signs}
We orient all moduli spaces $\mathcal M^{d+1}(y_0,\dots,y_d)$, $d \geq
1$, by choosing, as in Section \ref{sec:signs}, an orientation for
every $o(y)$, and requiring that the canonical identification
\begin{equation*}
    \Lambda^{\mathrm{top}}(T_{ u}\mathcal M^{d+1}( y_0,\dots, y_d)) \cong
    o( y_0)\otimes o( y_1)^\vee \otimes \dots \otimes o( y_d)^\vee
\end{equation*}
be orientation-preserving. In particular, this yields a sign
$\mathrm{sgn}(u) \in \{\pm 1\}$ for every regular curve $u \in
\mathcal M^{d+1}( y_0,\dots, y_d)$ of index zero.

To express the $A_\infty$-compositions of $\mathscr B_{\gamma,h}$
(partially) in terms of those of $\FukExt(M)$, we must orient the
moduli spaces $\mathcal M^{d+1}(\widetilde y_0,\dots,\widetilde y_d)$
in a compatible way. 
The requirement is as follows: Given some $\widetilde y =
(o_j,y)$ with odd $j$ and an orientation of $o(y)$ as above, we orient
$o(\widetilde y)$ such that the canonical identification
\begin{equation*}
    o(\widetilde y) \cong o(y)
\end{equation*}
from Lemma \ref{lem:comp-of-indices-and-orientations} matches up the
orientations. (For the $\widetilde y = (o_j,y)$ with even $j$ it will
not matter how the $o(\widetilde y)$ are oriented.) The orientations
of the $\mathcal M^{d+1}(\widetilde y_0,\dots,\widetilde y_d)$ and in
particular the signs $\mathrm{sgn}(\widetilde u) \in \{\pm 1\}$ of
isolated curves $\widetilde u$ are then determined as described above.

\begin{lem}
    \label{lem:signs}
    Given intersection points $y_0,\dots,y_d$ and corresponding
    $\widetilde y_0,\dots, \widetilde y_d$ with $\widetilde y_i =
    (o_j,y_i)$ for some \emph{odd} $j$, the canonical identification
    $\mathcal M^{d+1}(y_0\dots,y_d) \cong \mathcal M^{d+1}(\widetilde
    y_0,\dots,\widetilde y_d)$ is orien\-ta\-tion-pre\-serv\-ing. In
    particular, $\mathrm{sgn}(u) = \mathrm{sgn} (\widetilde u)$ for
    all isolated curves $u$ and $\widetilde u = (o_j,u)$.
\end{lem}
\begin{proof}
    The statement follows from the commutativity of the diagram
    \begin{equation*}
        \xymatrix{
            \Lambda^{\mathrm{top}}(T_{ u}\mathcal M^{d+1}( y_0,\dots,
            y_d)) \ar[r]^{\cong}\ar[d]^{\cong} & 
            \Lambda^{\mathrm{top}}(T_{\widetilde u}\mathcal M^{d+1}(\widetilde y_0,\dots,\widetilde y_d))\ar[d]^{\cong}\\
            \Lambda^{\mathrm{top}}(\mathrm{ker}\,D_u) \ar[r]^{\cong}\ar[d]^{\cong} & 
            \Lambda^{\mathrm{top}}(\mathrm{ker}\,D_{\widetilde u}\ar[d]^{\cong})\\
            o(y_0) \otimes o(y_1)^\vee \otimes \dots \otimes
            o(y_d)^\vee \ar[r]^{\cong} & 
            o(\widetilde y_0)\otimes o(\widetilde y_1)^\vee \otimes \dots \otimes
            o(\widetilde y_d)^\vee
        }
    \end{equation*}
    in which the vertical isomorphisms are the canonical ones, and
    where the first row is induced by $\mathcal M^{d+1}(y_0\dots,y_d)
    \cong \mathcal M^{d+1}(\widetilde y_0,\dots,\widetilde y_d)$, the
    second by the splitting $D_{\widetilde u} = D_{\widetilde
        u,\mathbb C} \oplus D_u$ and the fact that
    $\mathrm{ker}\,D_{\widetilde u,\mathbb C} = 0$, and the third by
    the canonical identifications from Lemma
    \ref{lem:comp-of-indices-and-orientations}. 
\end{proof}

\subsection{Definition of the functors $c_j:\mathscr
    B_{\gamma,h}\to\FukExt(M)$.}
\label{sec:a_infty-quasi-isom}
We now complete the definition of the (presumable) $A_\infty$-functor
\begin{equation*}
    c_j \equiv c_{\gamma,h,j}: \mathscr B_{\gamma,h} \to \FukExt(M)
\end{equation*}
for odd $j$ with $1 \leq j \leq \ell$. On objects, $c_j$ takes
$(\widetilde L^\deco,\widetilde E)$ to $ (L^\deco,E)$, as already
mentioned in Section \ref{sec:quasi-isom-mathscr}. As for morphisms,
consider objects $\widetilde L_i \equiv (\widetilde
L_0^\deco,\widetilde E_i)$ for $i = 0,1$ and the corresponding $L_i
\equiv (L_i,E_i)$. Note that for each of the summands $CF(\widetilde
L_0, \widetilde L_1)^{o_j} \subset CF(\widetilde L_0, \widetilde L_1)$
appearing in the splitting \eqref{eq:splitting-of-auxiliary-homs} we
have a canonical isomorphism of vector spaces
\begin{equation*}
    CF(\widetilde L_0, \widetilde L_1)^{o_j} \cong CF(L_0,L_1).
\end{equation*}
This is compatible with gradings, because the indices of corresponding
intersection points $y$ and $\widetilde y = (o_j,y)$ satisfy
$i(\widetilde y) = i(y)$ whenever $j$ is odd (cf. Lemma
\ref{lem:comp-of-indices-and-orientations}). We define the first order
component of $c_j$ on morphisms to be the composition
\begin{equation*}
    c_j^1: CF(\widetilde L_0^\deco, \widetilde L_1^\deco) \to
    CF(\widetilde L_0^\deco,\widetilde L_1^\deco)^{o_j}
    \cong CF(L_0^\deco,L_1^\deco),
\end{equation*}
where the first map is the projection onto morphisms of type
$o_j$. The higher order components $c_j^d$ for $d \geq 2$ are defined
to be identically zero.

\begin{prop}[Cf. {\cite[Prop.\ 4.2.3]{Biran-Cornea--Lag-Cob-II}}]
    \label{prop:relating-Fuk-and-B}
    The functors $c_j:\mathscr B_{\gamma,h} \to \FukExt(M)$ are
    $A_\infty$-quasi-isomorphisms. Moreover, $c_i$ and $c_j$ are
    homotopic for any two odd $i,j$ with $1 \leq i,j \leq \ell$ .
\end{prop}

\begin{proof}
    The proof is an adaptation of the proof of Proposition 4.2.3 in
    \cite{Biran-Cornea--Lag-Cob-II} that takes into account gradings
    and signs. We have already noted that the maps $c_j^1$ are
    compatible with gradings. Using the arguments in
    \cite{Biran-Cornea--Lag-Cob-II} together with Lemma
    \ref{lem:signs}, one sees that the differentials of $CF(\widetilde
    L_0^\deco,\widetilde L_1^\deco)$ and $CF(L_0^\deco,L_1^\deco)$ are
    related as described in \cite[Remark
    4.2.2]{Biran-Cornea--Lag-Cob-II}, which implies that the complexes
    are quasi-isomorphic. Then one checks that the $\{c_j^d\}$ really
    define $A_\infty$-functors $c_j$ with a common homotopy inverse $e
    \equiv e_{\gamma,h}: \FukExt(M) \to \mathscr B_{\gamma,h}$.
\end{proof}

\subsection{Exact triangles from cobordisms.}
\label{sec:constr-exact-triangl}
The construction of the exact triangles associated to a Lagrangian
cobordism now follows the scheme in \cite[Sections 4.3,
4.4]{Biran-Cornea--Lag-Cob-II}. Since even to outline this would
require quite a bit of additional notation, we content ourselves with
saying that all relevant statement carry over with minor
modifications. The bridge between $\FukExt(M)$ and the world of
cobordisms is now provided by the $A_\infty$-quasi-isomorphisms
\begin{equation*}
    \mathscr B_{\gamma,h}~
    \raisebox{-3pt}{$\stackrel{\xrightarrow{~~c_j~~}}{\xleftarrow[~~\,e\,~~]{}}
        $}
    ~\FukExt(M)
\end{equation*}
from Proposition \ref{prop:relating-Fuk-and-B}. To state the upshot of
all this, let $L^\deco \equiv (L^\deco,E)$ and $L_1^\deco \equiv
(L_1^\deco,E_1), \dots, L_k^\deco \equiv (L_k^\deco,E_k)$ be
Lagrangian branes in $M$ carrying local systems, and let
\begin{equation*}
    V^\deco: L^\deco \leadsto (L_1^\deco,\dots,L_k^\deco)
\end{equation*}
be a Lagrangian cobordism equipped with a brane structure and a local
system which restrict to the given ones on the ends.

\begin{prop}
    \label{prop:cob-cone-decomp}
    There exist $\FukExt(M)$-modules $\mathcal M_{V^\deco}^1,\dots,
    \mathcal M_{V^\deco}^k$ such that $\mathcal M_{V^\deco}^1 =
    \mathrm{Yon}(L_1^\deco)$ and such that there are exact triangles
    \begin{equation*}
        \mathcal M_{V^\deco}^{j-1} \to \mathrm{Yon}(L_j^\deco) \to \mathcal
        M_{V^\deco}^j \to \mathcal M_{V^\deco}^{j-1} [1]
    \end{equation*}
    in $\mathrm{mod}(\FukExt(M))$ for $j = 2,\dots,k$. Moreover,
    there is an $A_\infty$-quasi-isomorphism $\mathrm{Yon}(L^\deco)
    \to \mathcal M_{V^\deco}^k.$
\end{prop}

Here $\mathrm{Yon}: \FukExt(M) \to \mathrm{mod}(\FukExt(M))$ denotes
the Yoneda embedding of $\FukExt(M)$ into the $A_\infty$-category
$\mathrm{mod}(\FukExt(M))$ of $A_\infty$-modules over itself. The
Proposition summarizes what would be the analogues of Corollary 4.3.3
and Proposition 4.4.1 of \cite{Biran-Cornea--Lag-Cob-II} in our
setting (and the $\mathcal M_{V^\deco}^j$ are the analogues of the
$\mathcal M_{V,\gamma_j,h_j}$ there). Note that the directions of the
arrows in the exact triangles are reversed compared to
\cite{Biran-Cornea--Lag-Cob-II} because we use cohomological
conventions.

Theorem \ref{thm:cob-cone-decomp} is an immediate consequence of
Proposition \ref{prop:cob-cone-decomp}.

\subsection{Cobordism and Grothendieck groups.}
\label{sec:cobordism-group-k}
Denote by $\Fuk(M)$ the full subcategory of $\FukExt(M)$ consisting of
Lagrangian branes with trivial rank one local systems. Let $\cobGr(M)$
and $\cobGrExt(M)$ be the Lagrangian cobordism group defined as in
Sections \ref{sec:lagr-cobord-group-1} and
\ref{sec:relations-from-local-systems} with respect to $\Lags =
\mathrm{Ob}\, \Fuk(M)$, $\LagsExt = \mathrm{Ob}\,\FukExt(M)$ and
relation subgroup $R$, $R^\sesLS$ induced by appropriate Lagrangians
cobordisms between these generators. Denote by $K_0(D \Fuk(M))$ and
$K_0(D \FukExt(M))$ be the Grothendieck groups of the derived Fukaya
categories.
\begin{prop}
    \label{prop:Cob-K-group-hom} There exist canonical surjective
    group homomorphisms
    \begin{equation*}
        \Theta: \cobGr(M) \to K_0(D \Fuk(M)) ~ \text{and} ~
        \Theta^\deco: \cobGrExt(M) \to K_0(D \FukExt(M))
    \end{equation*}
    induced by $L \mapsto L$ resp. $(L,E) \mapsto (L,E)$.
\end{prop}

The statement for $\cobGr(M) \to K_0(D \Fuk(M))$ follows immediately
from Theorem \ref{thm:cob-cone-decomp} (together with Lemma
\ref{lem:linearization-and-K0}), while that for $\cobGrExt(M) \to
K_0(D \FukExt(M))$ needs in addition Proposition
\ref{prop:ses-ls-exact-triangle}, which implies that the relations in
$\cobGrExt(M)$ coming from short exact sequences of local systems are
respected.

\section{Lagrangian surgery and cobordisms in $\widetilde T^2$}
\label{sec:lagr-cobord-T2}
Surgering transversely intersecting Lagrangians $L_0,L_1$ means
cutting out small neighbourhoods in $L_0$ and $L_1$ of each
intersection point and gluing in Lagrangian handles, such as to
produce a new Lagrangian submanifold $L_0 \# L_1$ (see
\cite{Polterovich--Surgery-of-Lagrangians}). The local model is the
surgery of the Lagrangian subspaces $\mathbb R^n$ and $i \mathbb R^n$
of $\mathbb C^n$. Biran--Cornea \cite{Biran-Cornea--Lag-Cob-I}
describe how to construct a Lagrangian cobordism
\begin{equation*}
    \mathbb R^n \# i \mathbb R^n \leadsto (\mathbb R^n, i \mathbb
    R^n).
\end{equation*}
In general, one can produce cobordisms $L_0 \# L_1 \leadsto (L_0,L_1)$
by gluing in this local model via Darboux charts. We describe this
construction for $n = 1$ in Sections \ref{sec:case-n=1} and
\ref{sec:surg-surf-cobord}. Then we address some specifics of the
cobordisms coming from surgery of curves in $T^2$. 

Throughout this section, we identify $\mathbb R^2 \cong \mathbb C$ in
the usual way, so that Lagrangian cobordisms now live in $\widetilde M
= \mathbb C \times M$.

\subsection{Model cobordisms in dimension 1.}
\label{sec:case-n=1}
\begin{figure}[t]
    \centering
    \includegraphics[scale=0.7]{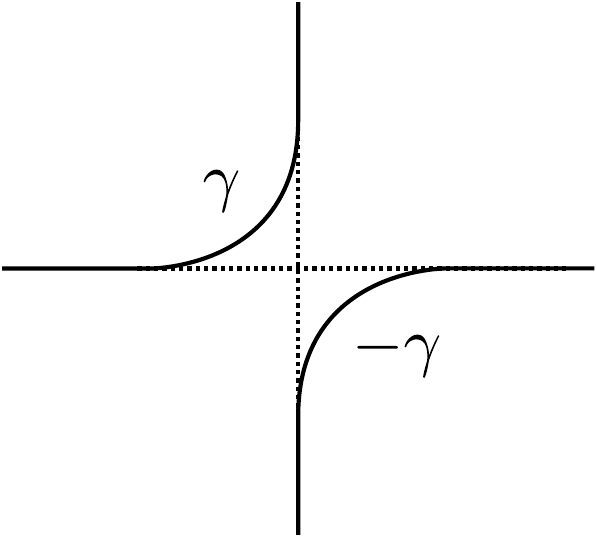}
    \caption{$\mathbb R \#_\gamma i\mathbb R$.}
    \label{fig:surgery}
\end{figure}

To surger $\mathbb R$ and $i \mathbb R \subset \mathbb C$, cut out two
neighbourhoods of $0 \in \mathbb C$ and connect the ends on $\mathbb
R$ thus created with the ends on $i \mathbb R$ by two curve
segments. More formally, choose a smooth curve $\gamma: \mathbb R \to
\mathbb C$, $\gamma = a + ib$, such that 
\begin{itemize}
    \label{conditions-on-gamma}
\item $\gamma(t) = t$ for $t \in (-\infty,-\varepsilon]$,
\item $\gamma(t) = it$ for $t \in [\varepsilon,\infty)$,
\item $a'(t),\, b'(t) > 0$ for $t \in (-\varepsilon, \varepsilon)$
\end{itemize}
for some $\varepsilon > 0$. Then define $\mathbb R \#_\gamma i \mathbb
R = \{ \pm \gamma(t) ~|~ t \in \mathbb R\}$, as shown in Figure
\ref{fig:surgery}.

To construct the corresponding Lagrangian cobordism
\begin{equation*}
    V_\gamma: \mathbb R \#_\gamma i\mathbb R \leadsto (\mathbb R, i\mathbb R),
\end{equation*}
consider the embedding $\psi_\gamma: \mathbb R \times S^1 \to\mathbb
C^2$, $(t, (x_1,x_2)) \mapsto (\gamma(t) \cdot x_1, \gamma(t) \cdot
x_2)$, where $S^1 = \{(x_1,x_2) \in \mathbb R^2 ~|~ x_1^2 + x_2^2 =
1\}$. As a first step, we set
\begin{equation*}
    V_\gamma' = \psi_\gamma (U)
\end{equation*}
for $U = \big\{(t, (x_1,x_2)) \in \mathbb R \times S^1 ~|~ 0 \leq x_1,
-2 \varepsilon \leq tx_1 \leq 2 \varepsilon \big\} \subset \mathbb R
\times S^1.$ 
By checking where $\psi_\gamma$ takes the boundary components of $U$,
one sees that $V_\gamma'$ is a manifold with boundary
\begin{equation*}
    \partial V_\gamma' = \{-2 \varepsilon\} \times
    \mathbb R \,\cup\, \{2 \varepsilon i\} \times i \mathbb R \,\cup\, \{0\}
    \times \mathbb R \#_\gamma i\mathbb R \,\subset\, \mathbb C \times
    \mathbb C.
\end{equation*}
To complete the construction of $V_\gamma$, we extend the part of
$\partial V_\gamma'$ lying over $0 \in \mathbb C$ to a cylindrical end
as explained in the proof of Lemma 6.1.1 of
\cite{Biran-Cornea--Lag-Cob-I} (the other two ends are already
cylindrical).
\begin{figure}[t]
    \centering
    \includegraphics[scale=1]{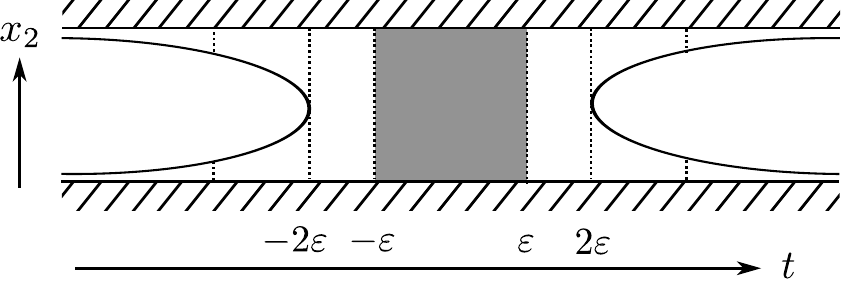}
    \caption{$V_\gamma$ schematically.}
    \label{fig:model-cob}
\end{figure}

Figure \ref{fig:model-cob} shows a schematic picture of $V_\gamma$
that should be viewed in light of the parametrisation $\psi_\gamma: U
\to V_\gamma'$. To understand the labels, note that the $U$ can be
identified with a subset of $\mathbb R \times [-1,1]$ via
$(t,(x_1,x_2)) \mapsto (t,x_2)$. The striped regions represent the
cylindrical end we attached at the very end of the construction, while
the region between them represents $V_\gamma'$. The two curved lines
correspond to the negative ends $\mathbb R$ and $i\mathbb R$ of the
cobordism, and the horizontal ones to the two components of its
positive end $\mathbb R \# i\mathbb R$.

\begin{figure}[t]
    \centering
    \includegraphics[scale=0.7]{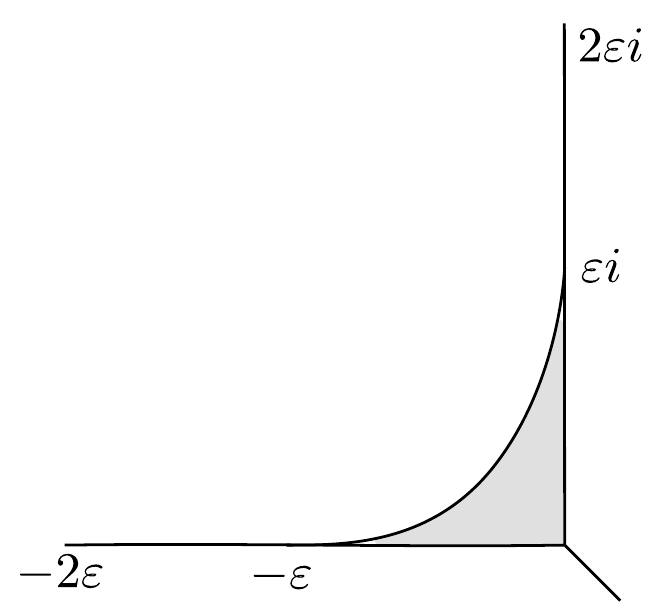}
    \caption{The projection of $V_\gamma$ to $\mathbb C$.}
    \label{fig:model-cob-proj}
\end{figure}

Figure \ref{fig:model-cob-proj} shows the projection of $V_\gamma$
under $\pi : \mathbb C^2 \to \mathbb C$ (the projection to the first
factor). In this picture, the line segment in the lower right is the
image of the cylindrical end attached at the very end, the rest is the
image of
$V_\gamma'$. 

We will refer to the following subsets of $V_\gamma$ as its
\textit{fat} resp. its \textit{thin} part:
\begin{equation*}
    \begin{aligned}
        V_\gamma^f =&~ \psi_\gamma \big(\{(t,(x_1,x_2)) \in \mathbb R
        \times
        S^1~|~ 0 \leq x_1, -\varepsilon \leq t \leq \varepsilon \}\big),\\
        V_\gamma^t =&~ \psi_\gamma\big( \{(t,(x_1,x_2)) \in \mathbb R
        \times S^1 ~|~ x_1 = 0, \,x_2 = \pm1 \} \big). 
    \end{aligned}
\end{equation*}
In Figure \ref{fig:model-cob}, $V_\gamma^f$ is represented by the
shaded square, and $V_\gamma^t$ is represented by the union of the two
horizontal lines. Note that $V_\gamma^f$ is precisely the part of
$V_\gamma$ lying over the ``fat'' part of $\pi(V_\gamma)$ (the shaded
region in Figure \ref{fig:model-cob-proj}), whereas $V_\gamma^t$
projects via $\pi$ to $0 \in \mathbb C$. Moreover, note that
$V_\gamma$ deformation retracts onto $V_\gamma^f \cup V_\gamma^t$.

\subsection{Surgery in surfaces and cobordisms.}
\label{sec:surg-surf-cobord}
Let now $L_0$ and $L_1$ be two (possibly disconnected) Lagrangians in
a surface $\Sigma$ intersecting transversely in $m$ points
$q_0,\dots,q_{m-1}$. Fix small neighbourhoods $U_j \subset \Sigma$ of
the $q_j$ and Darboux charts $\phi_j: U_j \to \mathbb C$ such that
$\phi_j (L_0 \cap U_j) \subset \mathbb R$ and $\phi_j(L_1 \cap U_j)
\subset i \mathbb R$. Moreover, choose curves $\gamma_j: \mathbb R \to
\mathbb C$ as in the previous section. 
The surgered manifold $L_0 \# L_1$ corresponding to these data is
obtained by gluing in the local surgery models corresponding to the
$\gamma_j$ using the charts $\phi_j$.

To construct the corresponding cobordism $V:L_0 \# L_1 \leadsto
(L_0,L_1)$, we first define subsets of $\mathbb C$ as follows:
\begin{equation*}
    \begin{aligned}
        I_0 &= \{x\in \mathbb R~|~ -2\varepsilon \leq x \leq 0\},\\
        I_1 &= \{iy \in i\mathbb R~|~ 0 \leq y \leq 2 \varepsilon\},\\
        I_2 &= \{ x-ix \in \mathbb C~|~ 0 \leq x \leq 2 \varepsilon\}.
    \end{aligned}
\end{equation*}
Moreover, we set $\tilde L_i = L_i \setminus \bigcup_{j=0}^{m-1}(L_i
\cap \phi_j^{-1}(B_{2 \varepsilon}(0)))$ for $i = 0,1$, and $\tilde
V_{\gamma_j} = V_{\gamma_j} \cap (\mathbb C \times B_{2
    \varepsilon}(0))$ for $j = 0,\dots,m-1$, where the $V_{\gamma_j}$
are the model cobordisms from the previous section. Then we set
\begin{equation*}
    V = (I_0 \times \tilde L_0) ~ \cup ~ (I_1 \times  \tilde L_1)
    ~\cup~  (I_2 \times (\tilde L_0 \cup \tilde L_1 )) ~\cup ~
    \bigcup_{j=0}^{m-1} (\mathrm{id} \times \phi_j^{-1})(\tilde V_{\gamma_j}),
\end{equation*}

\begin{rmk}
    \label{rmk:surgery-and-orientations}
    For this to work, it was necessary to choose \emph{all} charts
    $\phi_j$ such that $L_0 \cap U_j$ was mapped to $\mathbb R$ and
    $L_1 \cap U_j$ to $i \mathbb R$. This reflects that one needs to
    select one of the $L_i$ whose corresponding cylindrical end lies
    over $\mathbb R$, and the other for which it lies over $i \mathbb
    R$. 
    Of course, one can swap $\mathbb R$ and $i\mathbb R$ by a
    $\frac{\pi}{2}$-rotation, so that this requirement does not impose
    any restriction at first. However, if one wants to glue in all
    local surgery models in a way compatible with given orientations
    of the $L_i$ (which implies that the resulting cobordism is
    orientable), then the local topological type of the surgery at one
    intersection point determines that at all others.
    \begin{figure}
        \captionsetup{width=0.75\textwidth}
        \centering
        \includegraphics[scale=0.7]{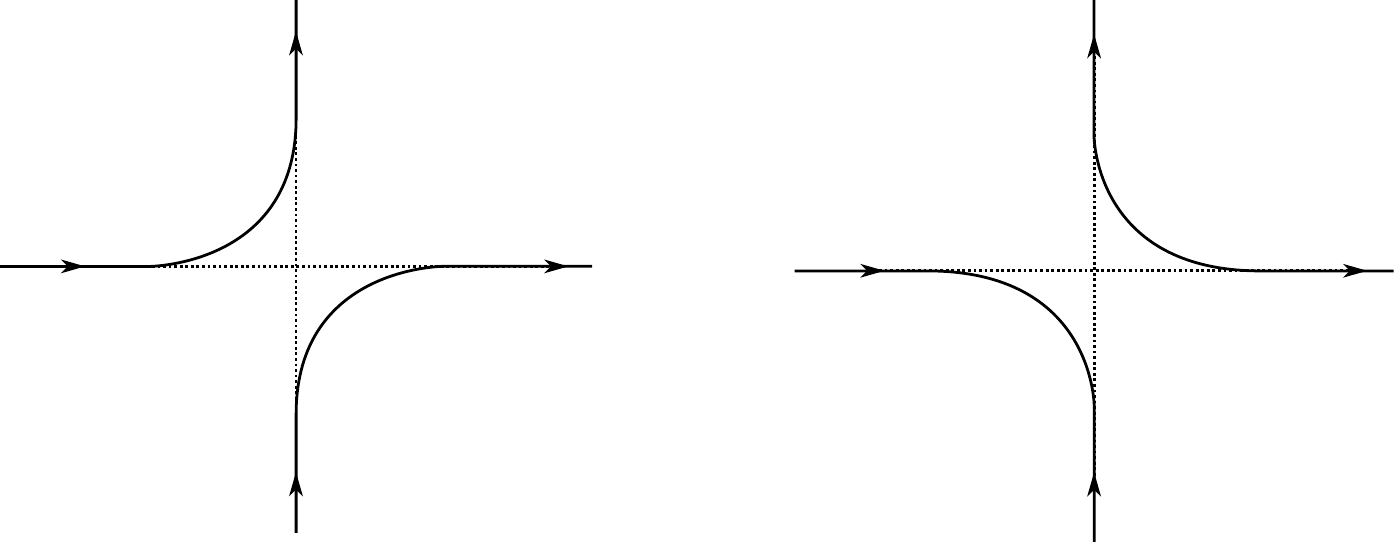}
        \caption{Surgering compatibly and non-compatibly with given
            orientations.}
        \label{fig:surgery-local-types}
    \end{figure}
\end{rmk}

As in the local case, we define \textit{fat} and \textit{thin} parts
of $V$ by
\begin{equation*}
    \begin{aligned}
        V^f =&~ \bigcup_{j=0}^{m-1} (\mathrm{id} \times
        \phi_j^{-1})(\tilde V_{\gamma_j}^f),\\
        V^t =&~ \bigcup_{j=0}^{m-1} (\mathrm{id} \times
        \phi_j^{-1})(\tilde V_{\gamma_j}^t) ~\cup~ \{0\} \times
        (\tilde L_0 \cup \tilde L_1),
    \end{aligned}
\end{equation*}
where $\tilde V_{\gamma_j}^f = V_{\gamma_j}^f \cap (\mathbb C \times
B_{2\varepsilon}(0))$ and $\tilde V_{\gamma_j}^t = V_{\gamma_j}^t \cap
(\mathbb C \times B_{2\varepsilon}(0))$. One can deformation retract
$V$ onto $V^f \cup V^t$ by retracting the glued-in $\tilde
V_{\gamma_j}$ onto their fat and thin parts, and the remaining
cylindrical parts of $V$ onto $\{0\} \times (\tilde L_0 \cup \tilde
L_1)$. From this description one sees readily that $V$ is homotopy
equivalent to a 1-dimensional CW-complex with $m$ cells of dimension 0
and $2m$ cells of dimension 1, where $m = \# \,L_0 \cap L_1$.

\subsection{Cobordisms in $\widetilde T^2$: Maslov class and
    grading.}
\label{sec:maslov-class}
We identify the tangent bundle of $\widetilde T^2 = \mathbb C \times
T^2$ with $\widetilde T^2 \times \mathbb R^4$ in the obvious way, so
that the Gau\ss~ map of a Lagrangian in $\widetilde T^2$ takes values
in the Lagrangian Grassmannian $\mathrm{Gr}(\mathbb R^4)$.

In the following, we call a Lagrangian $L \subset T^2$ \emph{straight}
if it lifts to a straight line in the universal cover. The cobordisms
of relevance for us come from surgering oriented straight Lagrangians
in $T^2$ in a way that is compatible with the orientations in the
sense of Remark \ref{rmk:surgery-and-orientations}. (This ensures not
only orientability of the resulting cobordisms, but also that $L_0
\#L_1$ has no contractible components.)

\begin{lem}
    \label{lem:Gauss-cobsT2}
    Let $L_0, L_1$ be oriented straight Lagrangians in $T^2$ of
    different slopes, let $L_0\#L_1$ be the result of an
    orientation-compatible surgery, and let $V: L_0 \# L_1 \leadsto
    (L_0,L_1)$ be the resulting cobordism. Then the Gau\ss~ map
    $\Gamma_V: V \to \mathrm{Gr}(\mathbb R^4)$ is null-homotopic. In
    particular, $V$ has vanishing Maslov class.
\end{lem} 

\begin{proof} We will show that $\Gamma_V$ is homotopic to a map which
    factors through a map defined on a contractible domain; this
    statement immediately implies the assertion. Note that we need
    only show this for $\Gamma_V \vert_{V^f \cup V^t}$, since $V$
    deformation retracts onto $V^f \cup V^t$. Moreover, we can assume
    that the local surgeries at all intersection points $q_i$ are
    defined with respect to the same local model, and thus that $V$ is
    obtained by gluing in the same $V_\gamma$ at every $q_i$; namely,
    isotopies of the curves $\gamma_j$ (with respect to which the
    surgery is defined) to one given curve $\gamma$ induce a
    Lagrangian isotopy of the resulting cobordisms.
    
    Under these assumptions, all maps $\Gamma_V \circ (\mathrm{id}
    \times \phi_j^{-1}): V_\gamma^f \to \mathrm{Gr}(\mathbb R^4)$ are
    equal to one and the same map $\tilde \Gamma:V_\gamma^f \to
    \mathrm{Gr}(\mathbb R^4)$. Since $\Gamma_V$ is constant on every
    connected component of $(V^f \cup V^t) \setminus V^f$, it factors
    through $\tilde \Gamma$, whose domain $V_\gamma^f$ is
    contractible.
\end{proof}

A grading of a Lagrangian $L \subset \widetilde T^2$ is a function
$\alpha: L \to \mathbb R$ lifting the composition $\mathrm{det}^2
\circ \Gamma_L: L \to \mathrm{Gr}(\mathbb R^4) \to S^1$ (cf. Section
\ref{sec:gradings}), where $\Gamma_L: L \to \mathrm{Gr}(\mathbb R^4)$
denotes the Gau\ss~ map. Note that every Lagrangian cobordism $V: L_0
\# L_1 \leadsto (L_0,L_1)$ coming from surgery of two linear
Lagrangians in $T^2$ admits a grading $\alpha_V: V \to \mathbb R$
because its $\Gamma_V$ is null-homotopic by Lemma
\ref{lem:Gauss-cobsT2}. The next lemma tells how the restrictions of
$\alpha_V$ to the ends are related.

\begin{lem}
    \label{lem:gradings-cobsT2}
    The restrictions of the grading $\alpha_V: V \to \mathbb R$ to the
    ends satisfy $\alpha_V \vert_{L_0} < \alpha_V \vert_{L_0 \# L_1} <
    \alpha_V \vert_{L_1}$.
\end{lem}

Here we assume the positive end $L_0 \# L_1$ has been ``linearized''
by a Hamiltonian isotopy, so that $\alpha_V\vert_{L_0\#L_1}$ is
constant. The lemma can be proven easily by examining the local models
from Section \ref{sec:case-n=1}.

\section{Preliminaries on $\cobGr(T^2)$}
\label{sec:comp-cobGr}
The Lagrangian cobordism group $\cobGr(T^2)$ we study is defined as in
Section \ref{sec:lagr-cobord-group-1} with $\mathscr L (T^2) =
\mathrm{Ob}\,\Fuk(T^2)$. That is, $\cobGr(T^2)$ has as generators
non-contractible curves equipped with brane structures, and relations
coming from cobordisms of vanishing Maslov class equipped with
compatible brane structures. In this section we prove as much as
possible about $\cobGr(T^2)$ as we can at this point. The result,
Proposition \ref{prop:an-exact-sequence}, will later be upgraded to
Theorem \ref{thm:MainThm-SES}.

\subsection{Notation for curves.}
\label{sec:notation}
\begin{figure}[t]
    \centering
    \includegraphics[scale=0.8]{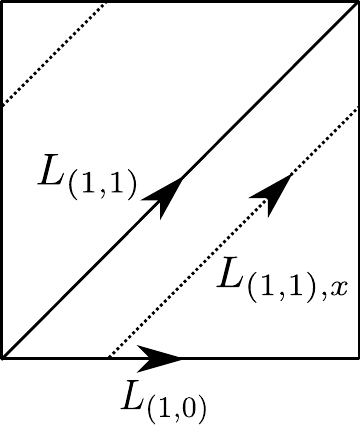}
    \caption{}
    \label{fig:Lags-in-T2}
\end{figure}

Let $m,n \in \mathbb Z$ such that $\mathrm{gcd}(m,n) = 1$, and let $x
\in \mathbb R$. If $(m,n) \neq (\pm 1,0)$, we define
\begin{equation*}
    L_{(m,n),x} \subset T^2
\end{equation*}
to be the oriented straight curve of slope $(m,n) \in \mathbb Z^2
\cong H_1(T^2;\mathbb Z)$ passing through the point $(x,0) \in
T^2$. We define $L_{(\pm 1,0),x}$ to be the straight horizontal curve
through $(0,x) \in T^2$ oriented such that it represents $(\pm 1,0)
\in H_1(T^2;\mathbb Z)$. When $x = 0$, we abbreviate
\begin{equation*}
    L_{(m,n)} \equiv L_{(m,n),0}.
\end{equation*}
See Figure \ref{fig:Lags-in-T2} for an illustration. We view these
curves as objects of $\Fuk(T^2)$ by equipping them with the standard
grading and the bounding $Pin$ structure, which we refer to as the
\emph{standard brane structure}; note that the orientation induced on
$L_{(m,n)}$ by the standard grading agrees with the orientation it had
a priori (cf. Section \ref{sec:fukaya-category}).

\subsection{Cylinders.}
\label{sec:cylinders}
A cylinder in $T^2$ is a smooth map $u: C = [0,1] \times S^1\to
T^2$. Given such a cylinder $u$, we equip its boundary curves $L^i =
u(\{i\} \times S^1) \subset T^2$, $i = 0,1$, with the boundary
orientation induced from the standard orientation of $[0,1] \times
S^1$, and with the standard brane structure. Note that the area
$\int_Cu^*\omega \in \mathbb R$ of $u$ is determined up to an integer
by its boundary, as can be seen by considering the cover of $T^2$
corresponding to the subgroup of $\pi_1(T^2)$ generated by $u$. This
justifies the wording in the following lemma.

\begin{lem}
    \label{lem:can-group-hom}
    There exists a canonical group homomorphism $\zeta: \mathbb
    R/\mathbb Z \to \cobGr(T^2)$ taking $x \in \mathbb R/\mathbb Z$ to
    $[L^1] - [L^0] \in \cobGr(T^2)$, where $L^0, L^1$ are any two
    non-contractible curves that bound a cylinder of area $x$.
\end{lem}

\begin{proof} Our task is to prove that the map described is
    well-defined (i.e., independence of the chosen cylinder), and that
    it is a group homomorphism. Note that we may always replace curves
    by the straight representatives of their Hamiltonian isotopy
    classes, as Hamiltonian isotopies change neither areas nor
    cobordism classes.

    Our first claim is that the cobordism class of the boundary of a
    cylinder of area $x \in \mathbb R/\mathbb Z$ with boundary circles
    representing the class $(0,\pm 1) \in H_1(T^2;\mathbb Z)$ is
    well-defined, i.e., that in $\cobGr(T^2)$ we have
    \begin{equation}
        \label{eq:zeta-well-defd}
        [L_{(0,1),x + y}] - [L_{(0,1),y}] = [L_{(0,1),x}] - [L_{(0,1)}]
    \end{equation}
    for all $x,y \in \mathbb R/\mathbb Z$. Before we justify that,
    note that subtracting $[L_{(0,-1)}]$ from both sides of
    \eqref{eq:zeta-well-defd} and rearranging yields
    \begin{equation}
        \label{eq:zeta-group-hom}
        \Big([L_{(0,1),x + y }] - [L_{(0,1)}]\Big) = 
        \Big([L_{(0,1),x}] - [L_{(0,1)}]\Big) +
        \Big([L_{(0,1),y}] - [L_{(0,1)}]\Big),
    \end{equation}
    which shows that our map (if well-defined) is a group
    homomorphism. 

    To obtain \eqref{eq:zeta-well-defd}, observe that surgering
    $L_{(1,0)}$ and $L_{(0,1),x+y}$ has the same result as surgering
    $L_{(1,0),-x}$ and $L_{(0,1),y}$ up to Hamiltonian isotopy (namely
    $L_{(1,1),\frac{1}{2} + x + y}$ in both cases, which is a right
    shift by $x+y$ of what one would get when surgering $L_{(1,0)}$
    and $L_{(0,1)}$). This observation tells us that $[L_{(1,0)}] +
    [L_{(0,1),x+y}] = [L_{(1,0),-x}] + [L_{(0,1),y}]$, or equivalently
    that
    \begin{equation}
        \label{eq:zeta-well-defd-2}
        [L_{(0,1),x+y}] - [L_{(0,1),y}] = [L_{(1,0),-x}] -
        [L_{(1,0)}]
    \end{equation}
    for all $x, y \in \mathbb R/\mathbb Z$. The identity
    \eqref{eq:zeta-well-defd} follows from that as the right-hand side
    of \eqref{eq:zeta-well-defd-2} is independent of $y$.
  
    To complete the proof of well-definedness, we have to show that
    for any two general curves $L^0$ and $L^1$ bounding a cylinder of
    area $x$ (after reversing the orientation on $L^0$), we have
    $[L^1] - [L^0] = [L_{(0,1),x}] - [L_{(0,1)}]$.  One way of seeing
    this is by thinking about how to build the $L^i$ using iterated
    surgery of copies of $L_{(1,0)}$ and translated copies of
    $L_{(0,1)}$. For example, if $L^0$ can be constructed by
    iteratively surgering $m$ copies of $L_{(1,0)}$ and $n$ copies of
    $L_{(0,1)}$, then $L^1$ can be constructed by iteratively
    surgering $m$ copies of $L_{(1,0)}$ and $n$ copies of
    $L_{(0,1),\frac{x}{n}}$.
    From this we obtain
    \begin{equation*}
        [L^1] - [L^0] = n\Big([L_{(0,1),\frac{x}{n}}] - [L_{(0,1)}] \Big) =
        [L_{(0,1),x}] - [L_{(0,1)}],
    \end{equation*}
    where the second equality follows from \eqref{eq:zeta-group-hom}.
    (Alternatively, one can directly deduce from
    \eqref{eq:zeta-well-defd} and \eqref{eq:zeta-well-defd-2} that the
    boundaries of all ``vertical'' and all ``horizontal'' cylinders of
    area $x$ represent the same class, and then get the general case
    by applying suitable maps $A \in \mathrm{SL}(2,\mathbb Z)$.)
\end{proof}

\subsection{An exact sequence.}
\label{sec:an-exact-sequence}
Denote by $\zeta: \mathbb R/\mathbb Z \to \cobGr(T^2)$ the map
described in the previous subsection, and by $\eta: \cobGr(T^2) \to
H_1(T^2;\mathbb Z)$ the canonical map given by $[L]_{\Omega} \mapsto
[L]_{H_1}$. The following proposition collects what we can say about
$\cobGr(T^2)$ at this point.

\begin{prop}
    \label{prop:an-exact-sequence}
    The sequence of group homorphisms
    \begin{equation*}
        \mathbb R/\mathbb Z \to \cobGr(T^2) \to H_1(T^2;\mathbb Z) \to 0
    \end{equation*}
    is exact.
\end{prop}
\begin{proof}
    It is clear that the canonical map $\eta: \cobGr(T^2) \to
    H_1(T^2;\mathbb Z)$ is surjective and that $\mathrm{im}\, \zeta
    \subseteq \mathrm{ker}\,\eta$. To prove $\mathrm{ker}\, \eta
    \subseteq \mathrm{im}\,\zeta$, consider the map 
        \begin{equation}
            \label{eq:map-comp-cobGr}
            \mathbb R/\mathbb Z \oplus H_1(T^2;\mathbb Z) \to \cobGr(T^2),
            (x,(m,n)) \mapsto \zeta(x) 
            + m[L_{(1,0)}] + n [L_{(0,1)}],
        \end{equation}
     where we identify $H_1(T^2;\mathbb Z) \cong \mathbb Z^2$ in the
    obvious way. Observe that $\cobGr(T^2)$ is generated by
    $[L_{(1,0)}]$ and the elements of the family $[L_{(0,1),x}]_{x \in
        \mathbb R/\mathbb Z}$, because every Lagrangian can be
    obtained by iteratively surgering these. An alternative set of
    generators of $\cobGr(T^2)$ is given by $[L_{(1,0)}], [L_{(0,1)}]$
    and the family $\zeta(x) = [L_{(0,1),x}] - [L_{(0,1)}]$, $x \in
    \mathbb R/\mathbb Z$, which shows that the map
    \eqref{eq:map-comp-cobGr} is surjective.
    Together with the fact that the composition $H_1(T^2;\mathbb Z)
    \hookrightarrow \mathbb R/\mathbb Z \oplus H_1(T^2;\mathbb Z) \to
    \cobGr(T^2)$ is a section of $\eta$, this implies that
    $\mathrm{ker}\,\eta \subseteq \mathrm{im}\,\zeta$.
\end{proof}

In order to upgrade Proposition \ref{prop:an-exact-sequence} to
Theorem \ref{thm:MainThm-SES}, it remains to prove the injectivity of
$\zeta: \mathbb R/\mathbb Z \to \cobGr(T^2)$. It is in fact easy to
rule out the existence of a Lagrangian cobordism $L \leadsto L'$ for
isotopic but not Hamiltonian isotopic curves: The existence of such a
cobordism would imply that $HF(N,L) \cong HF(N,L')$ for any other
curve $N$ (as a consequence of spelling out what Proposition
\ref{prop:cob-cone-decomp} says), contradicting that we have $HF(L,L)
\cong \mathbb Z^2$ and $HF(L,L') = 0$. However, there might be more
complicated relations leading to the identity $[L] - [L'] = 0 \in
\cobGr(T^2)$. The fact that we cannot rule these out directly is one
reason for the detour via homological mirror symmetry we will take in
the next section.

\section{Homological mirror symmetry for $T^2$}
\label{sec:homol-mirr-symm}

\subsection{Abouzaid-Smith's mirror functor.}
\label{sec:abouz-smiths-mirr}
Abouzaid-Smith prove in \cite{Abouzaid-Smith--HMS4torus} that the
split-closed derived Fukaya category $D^\pi \Fuk^\#(T^2)$ is
equivalent to the derived category $D^b(X)$ of the \emph{Tate curve}
$X$, which is an elliptic curve over the Novikov field $\Lambda$ given
by a specific Weierstrass equation (see
\cite{Abouzaid-Smith--HMS4torus}). 
To state the precise result, let $P_0 \in X$ be a base point, denote
by $\mathcal O(nP_0)$ the line bundle corresponding to the divisor $n
P_0$, and by $\mathcal O_{P_0}$ the skyscraper sheaf with
one-dimensional stalk supported at $P_0$.

\begin{thm}
    \label{thm:HMS-for-T2}
    There exists an equivalence of triangulated categories
    \begin{equation*}
        \Phi:
        D^b(X) \simeq D^\pi \FukExt(T^2)
    \end{equation*}
    which takes $\mathcal O(nP_0)$ to $L_{(1,-n)}$ for every $n \in
    \mathbb Z$, and $\mathcal O_{P_0}$ to $L_{(0,-1),\frac{1}{2}}$.
\end{thm} 

The proof is based on parts of Polishchuk-Zaslow's computations in
\cite{polishchuk-zaslow--HMS-elliptic-curve99}, who work over an
elliptic curve defined over $\mathbb C$. To make the connection to
the current setting, one studies $X$ by means of its analytification
\begin{equation*}
    X^{an} = \Lambda^* / q^{\mathbb Z},
\end{equation*}
the quotient of $\Lambda^* = \Lambda \setminus \{0\}$ by the discrete
subgroup generated by $q \in \Lambda$. $X^{an}$ can be given the
structure of a rigid-analytic space over $\Lambda$ (we refer the
reader to \cite{Fresnel-vanderPut--Rigid-anal-geom} for general
background on rigid-analytic geometry, and to
\cite{Fresnel-vanderPut--Rigid-anal-geom,
    Silverman--AdvTopicsArithmEllCurves} for specific information on
the Tate curve). The complex-analytic $\theta$-functions appearing in
\cite{polishchuk-zaslow--HMS-elliptic-curve99} can be interpreted as
formal power series and hence functions on $\Lambda^*$, which give
rise to sections of rigid-analytic vector bundles over $X^{an}$,
analogously to the complex case. A rigid-analytic GAGA principle 
says that the categories of coherent algebraic sheaves on $X$ and of
coherent analytic sheaves on $X^{an}$ are equivalent, which allows to
translate back to algebraic geometry. We will drop the notational
distinction between $X$ and $X^{an}$ in the following.

To outline Abouzaid-Smith's proof of Theorem \ref{thm:HMS-for-T2},
consider on the algebraic side the full subcategory $\Gamma \mathscr
A^\vee \subset D^b_\infty(X)$ consisting of the line bundles $\mathcal
O(nP_0)$, $n \in \mathbb Z$ (where $D^b_\infty(X)$ is a dg-enhancement
of $D^b(X)$, see \cite{Abouzaid-Smith--HMS4torus}); on the symplectic
side, consider the full $A_\infty$-subcategory $\Gamma \mathscr A
\subset \FukExt(T^2)$ with objects the Lagrangians $L_{(1,n)}$, $n \in
\mathbb Z$. Both collections of objects split-generate on their
respective sides. Polishchuk-Zaslow's computations
\cite{polishchuk-zaslow--HMS-elliptic-curve99} imply that the
cohomological categories $H^0(\Gamma \mathscr A^\vee)$ and $H^0(\Gamma
\mathscr A)$ are equivalent by a functor taking $\mathcal O(nP_0)$ to
$L_{(1,-n)}$. Then a deformation theoretic result by Polishchuk
\cite{Polishchuk--ExtOfHomCoordRings-2003} says that $H^0(\Gamma
\mathscr A^\vee) \simeq H^0(\Gamma \mathscr A)$ can be equipped with
an essentially unique non-formal $A_\infty$-structure. As the
$A_\infty$-structures on $\Gamma \mathscr A^\vee$ and $\Gamma \mathscr
A$ are both non-formal, Abouzaid-Smith conclude the existence of an
$A_\infty$-quasi-equivalence between $\Gamma \mathscr A^\vee$ and
$\Gamma \mathscr A$, which extends to an $A_\infty$-quasi-equivalence
between the respective split-closures. The triangulated equivalence of
Theorem \ref{thm:HMS-for-T2} is then obtained by taking cohomology.

\begin{rmk}
    Abouzaid-Smith \cite{Abouzaid-Smith--HMS4torus} really prove
    Theorem \ref{thm:HMS-for-T2} for $\Fuk(T^2)$, the version of the
    Fukaya category without local systems. However, the proof
    in \cite{Abouzaid-Smith--HMS4torus} goes through for both
    $\Fuk(T^2)$ and $\FukExt(T^2)$, because the $L_{(1,n)}$
    split-generate both versions. One can prove this, for example,
    with the machinery from Section \ref{sec:lagr-cobord-cones}, by
    iteratively surgering the $L_{(1,n)}$ and equipping the resulting
    cobordisms with suitable local systems. 
    Cf. Remark 6.7 in \cite{Abouzaid-Smith--HMS4torus}, which
    indicates that the effects of allowing appropriate non-trivial
    local systems and of taking the split-closure of $D \Fuk(T^2)$ are
    equivalent.
\end{rmk}

\subsection{Recovery of the mirror functor.}
\label{sec:recov-mirr-funct}
The triangulated equivalence $\Phi: D^b(X) \to D^\pi\FukExt(T^2)$ of
Theorem \ref{thm:HMS-for-T2} is established using a deformation
theoretic argument, and a priori it is not clear how precisely $\Phi$
acts on arbitrary objects. It is not even obvious that the mirror
object of every indecomposable sheaf is a Lagrangian with a local
system, as opposed to some ``abstract'' object introduced when passing
to the split-closure. In contrast, the equivalence constructed by
Polishchuk-Zaslow \cite{polishchuk-zaslow--HMS-elliptic-curve99} is
given by explicit formulae, but it is not clear how compatible it is
with the triangulated structure. 

The aim here is to partially recover the effect of $\Phi$ on
objects. 
We denote by $\mathcal V(r,d)$ the set of isomorphism classes of
indecomposable vector bundles of rank $r$ and degree $d$ on $X$, and
we identify $H_1(T^2;\mathbb Z) \cong \mathbb Z^2$ in the standard
way. As usual, we include objects of $\mathrm{Coh}\,X$ in $D^b(X)$ by
viewing them as complexes concentrated in degree 0.

\begin{prop}
    \label{prop:recov-mirr-funct}
    (i) Let $\mathscr S$ be an indecomposable skyscraper sheaf on $X$
    with stalk of rank $h$. Then $\Phi(\mathscr S)$ is isomorphic to a
    Lagrangian of slope $(0,-1) \in H_1(T^2;\mathbb Z)$ equipped with
    a local system of rank $h$.

    (ii) Let $E \in \mathcal V(r,d)$ and set $h=
    \mathrm{gcd}(r,d)$. Then $\Phi(E)$ is isomorphic to a Lagrangian
    of slope $\frac{1}{h} (r,-d) \in H_1(T^2;\mathbb Z)$, equipped
    with an indecomposable local system of rank $h$.
\end{prop}

Every object of $D^b(X)$ is isomorphic to a direct sum of shifted
copies of vector bundles and skyscraper sheaves (see, e.g., Corollary
3.15 in \cite{Huybrechts--FM-transforms-in-AG}). Since moreover the
shift functor of $D^\pi \FukExt(T^2)$ just shifts the grading of
Lagrangian branes, Theorem \ref{thm:HMS-for-T2} and Proposition
\ref{prop:recov-mirr-funct} together imply the following statement.

\begin{cor}
    \label{cor:objects-of-DF}
    Every object of $D^\pi \FukExt(T^2)$ is, up to isomorphism, a
    direct sum of objects of $\FukExt(T^2)$, i.e., of Lagrangian
    branes with local systems.
\end{cor}

\begin{cor}
    \label{cor:DF-is-equiv-to-DpiF}
    The inclusion $D \FukExt(T^2) \hookrightarrow D^\pi \FukExt(T^2)$
    is an equivalence.
\end{cor}

The proof of Proposition \ref{prop:recov-mirr-funct} will occupy the
rest of this section.

\subsection{Proof of Proposition \ref{prop:recov-mirr-funct} (i) for
    $h = 1$.} We start with some preliminary discussion. Every point
$Q \in X$ can be written as $Q = [-q^xM]$ with $x \in \mathbb R /
\mathbb Z$ and $M \in S^1\Lambda$ that are uniquely
determined.\footnote{Recall that by slight abuse of notation we write
    $X$ for the analytification $X^{an} = \Lambda^* / q^{\mathbb Z}$.}
From now on we fix the base point of $X$ to be
\begin{equation*}
    P_0 = [-q^{1/2}].
\end{equation*}
This choice determines an equivalence $\Phi: D^b(X) \to D^\pi
\FukExt(T^2)$ as described in Theorem \ref{thm:HMS-for-T2}, i.e., such
that $\Phi$ takes $\mathcal O(nP_0)$ to $L_{(1,-n)}$.

We denote by $O = [q^0]$ the neutral element for the natural group
structure on $X$ induced by multiplication on $\Lambda^*$. Since $X$
is an elliptic curve, it also carries an 
``elliptic curve group structure'' with of $O$ as neutral element, see
\cite[Section IV.4]{Hartshorne--AG}. In fact, both group structures
must coincide (cf.  \cite[p.127]{Fresnel-vanderPut--Rigid-anal-geom}),
and we denote the operation in this group by $\oplus$. It is related
to the operation in the divisor group by
\begin{equation}
    \label{group-op-elliptic-curve}
    Q \oplus Q' = Q'' \quad \Longleftrightarrow \quad Q + Q' \sim Q'' + O,
\end{equation}
where $\sim$ denotes linear equivalence of divisors (again, see
\cite[Section IV.4]{Hartshorne--AG}).

To determine the mirror images of the 1-dimensional skyscraper
sheaves, we will use that (almost) every such skyscraper sheaf can be
obtained as a direct summand of a cone on a morphism $\mathcal
O(-2P_0) \to \mathcal O_X$. To see this, note that
\begin{equation*}
    \mathrm{Hom}_{D^b(X)}(\mathcal O(-2P_0),\mathcal O_X) \cong
    H^0(X;\mathcal O(2 P_0)),
\end{equation*}
and hence the cone of every morphism $\mathcal O(-2P_0) \to \mathcal
O_X$ is of the form $\mathcal O_D$ for a divisor $D$ belonging to the
linear system $\vert 2 P_0\vert$. Observe that, by
\eqref{group-op-elliptic-curve}, $\vert 2 P_0 \vert$ is the set of all
divisors $D = Q + Q'$ such that $Q \oplus Q' = P_0 \oplus P_0$, hence
the set of divisors of the form
\begin{equation*} 
    D = [-q^xM] + [-q^{-x}M^{-1}]
\end{equation*}
with $x \in \mathbb R / \mathbb Z$ and $M \in S^1\Lambda$. Whenever
the points $[-q^xM]$ and $[-q^{-x}M^{-1}]$ are distinct, which is the
case unless $x \in \{ 0, \frac{1}{2}\}$ and $M \in \{ \pm q^0\}$, the
sheaf $\mathcal O_D$ corresponding to such $D$ is a direct sum of the
corresponding 1-dimensional skyscraper sheaves, that is,
\begin{equation*}
    \mathcal O_D = \mathcal O_{[-q^xM]} \oplus \mathcal O_{[-q^{-x}M^{-1}]}.
\end{equation*}
(And $\mathcal O_D = \mathcal O_{2[q^xM]}$ in the four cases in which
$[-q^xM] = [-q^{-x}M^{-1}]$.)

The space $\mathrm{Hom}_{D^b(X)}(\mathcal O(-2P_0),\mathcal O_X) \cong
H^0(X;\mathcal O(2P_0))$ has a preferred basis given by the theta
functions $\theta^0,\theta^1: \Lambda^* \to \Lambda$ defined by
\begin{eqnarray*}
    \theta^0 (w) &=& \sum_{n \in \mathbb Z} w^{2n}q^{n^2},\\
    \theta^1 (w) &=& \sum_{n \in \mathbb Z} w^{2n+1} q^{(n+\frac{1}{2})^2}.
\end{eqnarray*}
To understand this, observe that these functions satisfy the
functional equation $\theta^i(qw) = q^{-1}w^{-2}$, and therefore can
be considered as sections of the line bundle
\begin{equation*}
    \Lambda^* \times \Lambda / \big((w,\xi) \sim (qw, q^{-1}w^{-2} \xi)\big)
\end{equation*}
over $X$, which is $\mathcal O(2P_0)$. (In fact, the reason for
choosing the base point to be $P_0 = [-q^{1/2}]$ was to ensure that
$H^0(X;\mathcal O(2P_0))$ is spanned by these ``standard'' theta
functions.) We refer to \cite[Section
2.3]{polishchuk-zaslow--HMS-elliptic-curve99} or
\cite{Griffiths-Harris--Principles-of-AG} for the description of
vector bundles on elliptic curves (over $\mathbb C$) via
``multipliers'', to \cite[Section
4.7]{Fresnel-vanderPut--Rigid-anal-geom} for information on vector
bundles over rigid analytic spaces, and to \cite{Mumford--Tata-I} for
background on theta functions.

Now we turn to the symplectic side and study the cones on morphisms
\begin{equation*}
    c_1 \in \mathrm{Hom}_{D^\pi \FukExt(T^2)}(L_{(1,2)}, L_{(1,0)})
\end{equation*}
in $D^\pi\FukExt(T^2)$. This morphism space is spanned by the two
intersection points
\begin{equation*}
    c_1^0, c_1^1 \in L_{(1,2)} \cap L_{(1,0)},
\end{equation*}
see Figure \ref{fig:Cones-L(1,2)-L(1,0)} (more precisely, $c_1^0,
c_1^1$ are the identity homomorphisms between the fibres of the
respective local systems over these intersection points, which in our
model for local systems are all equal to $\Lambda$).

Recall that $L_{(1,2)}$ and $L_{(1,0)}$ are the mirror images of
$\mathcal O(-2P_0)$ and $\mathcal O_X$. Since the restriction of
$\Phi$ to the full subcategory of $D^b(X)$ consisting of the $\mathcal
O(nP_0)$ is essentially the Polishchuk-Zaslow functor
\cite{polishchuk-zaslow--HMS-elliptic-curve99}, the corresponding
isomorphism between morphism spaces is determined by
\begin{equation*}
    \mathrm{Hom}_{D^b(X)}(\mathcal O(-2P_0),\mathcal O_X) \ni
    \theta^i \leftrightarrow c_1^{i} \in \mathrm{Hom}_{D \FukExt(T^2)}(L_{(1,2)},L_{(1,0)})
\end{equation*}    
for $i = 0,1$.

\begin{lem}
    \label{lem:Cones-L(1,2)-L(1,0)}
    Let $c_1 = \sigma^0 c_1^0 + \sigma^1 c_1^1: L_{(1,2)} \to
    L_{(1,0)}$ be a non-zero morphism and let $s = \sigma^0 \theta^0 +
    \sigma^1 \theta^1$ be the corresponding section of $\mathcal
    O(2P_0)$. If $s$ vanishes at two distinct points $[-q^xM] \neq
    [-q^{-x}M^{-1}]$, then
    \begin{equation*}
        \mathrm{Cone}(c_1) \cong (L_{(0,-1),x}, E_M^1) \oplus (L_{(0,-1),-x},E_{M^{-1}}^1).
    \end{equation*}
    If $s$ has a double zero, then $\mathrm{Cone}(c_1) \cong
    (L_{(0,-1),x},E_M^2)$ for $x \in \{0,\frac{1}{2}\}$, $M \in \{\pm
    q^0\}$, and where $E_M^2$ is the unique non-trivial extension of
    $E_M^1$ by itself. 
\end{lem}

The lemma will be proven in the next subsection. Before that, we
finish proving Proposition \ref{prop:recov-mirr-funct}(i) for $h =
1$. Since $\Phi$ takes cones to cones, we can infer from Lemma
\ref{lem:Cones-L(1,2)-L(1,0)} and the discussion preceding it that
\begin{equation*}
    \Phi( \mathcal
    O_{[-q^xM]} \oplus \mathcal O_{[-q^{-x}M^{-1}]}) \cong (L_{(0,-1),x},E_M^1) \oplus (L_{(0,-1),-x},E_{M^{-1}}^1)
\end{equation*}
whenever $x \notin \{0,\frac{1}{2}\}$ and $M \notin \{\pm 1\}$. Since
there are no other ways of writing the object on the right-hand side
as a direct sum, we conclude that the mirror images of $\mathcal
O_{[-q^xM]}$ and $\mathcal O_{[-q^{-x}M^{-1}]}$ are
$(L_{(0,-1),x},E^1_M)$ and $(L_{(0,-1),-x},E^1_{M^{-1}})$. 
For the remaining four skyscrapers $\mathcal O_{[-q^xM]}$ with $x \in
\{0,\frac{1}{2}\}$ and $M\in \{\pm 1\}$, a similar argument shows that
they are mirror to the $(L_{(0,-1),x},E_M^1)$. Since every skyscraper
sheaf with stalk of rank $1$ is of the form $\mathscr S = \mathcal
O_{[-q^xM]}$, this concludes the proof.
\hfill $\square$

\subsection{Proof of Lemma \ref{lem:Cones-L(1,2)-L(1,0)}}
\label{sec:proof-lemma-on-cones} 
Set $Y_0 := L_{(1,2)}$ and $Y_1 := L_{(1,0)}$ and let $c_1 = \sigma^0
c_1^0 + \sigma^1 c_0^1 \in \mathrm{hom}^0(Y_0,Y_1)$ be a non-zero
morphism. We first rephrase what needs to be proven: Our task is to
show that there exists an object $Y_2$ as described in the statement
of the lemma, together with morphisms $c_2 \in
\mathrm{hom}^0(Y_1,Y_2)$ and $c_3 \in \mathrm{hom}^1(Y_2,Y_0)$, such
that the triangle
\begin{equation}
    \label{eq:surgery-exact-triangle}
    Y_0 \xrightarrow[]{[c_1]} Y_1 \xrightarrow[]{[c_2]} Y_2
    \xrightarrow[]{[c_3]} Y_0[1]
\end{equation}
is exact in $H(\FukExt(T^2))$. We will use Lemma 3.7 in
\cite{Seidel--Fukaya-Picard-Lefschetz-2008} to prove this, which says
that it suffices to show that
\begin{equation}
    \label{eq:exact-triangle-criterion}
    \begin{aligned}
        \mu^2(c_3,c_2) &\,=\, 0,\\
        \mu^2(c_1,c_3) &\,=\, 0,\\
        \mu^3(c_1,c_3,c_2) &\,=\, e_{Y_1},
    \end{aligned}
\end{equation}
where $e_{Y_1}$ is a chain representing the identity in
$\mathrm{Hom}^0(Y_1,Y_1)$, and to show the acyclity of the complex
$\big( \mathrm{hom}(X,Y_2)[1] \oplus \mathrm{hom}(X,Y_0)[1] \oplus
\mathrm{hom}(X,Y_1),\partial\big)$ for every test object $X$ (see
\cite{Seidel--Fukaya-Picard-Lefschetz-2008} for the description of the
differential $\partial$). We will in fact only verify parts of this
criterion, and then argue that this is already sufficient.

We first consider the case that $c_1$ is such that the corresponding
section $s = \sigma^0 \theta^0 + \sigma^1 \theta^1$ of $\mathcal O(2
P_0)$ vanishes at distinct points $[-q^{x}M],\,[-q^{-x}M^{-1}] \in X$,
as opposed to having a double zero. Set
\begin{equation}
    \label{eq:Y_2}
    Y_2 = (L_{(0,-1),x},E_M^1) \oplus (L_{(0,-1),-x},E_{M^{-1}}^1),
\end{equation}
and denote by $c_3^0, c_3^1$ the intersection points generating the
space $\mathrm{hom}^1(Y_2,Y_0)$, as depicted in Figure
\ref{fig:Cones-L(1,2)-L(1,0)}.
\begin{figure}[t]
    \centering
    \includegraphics[scale=0.65]{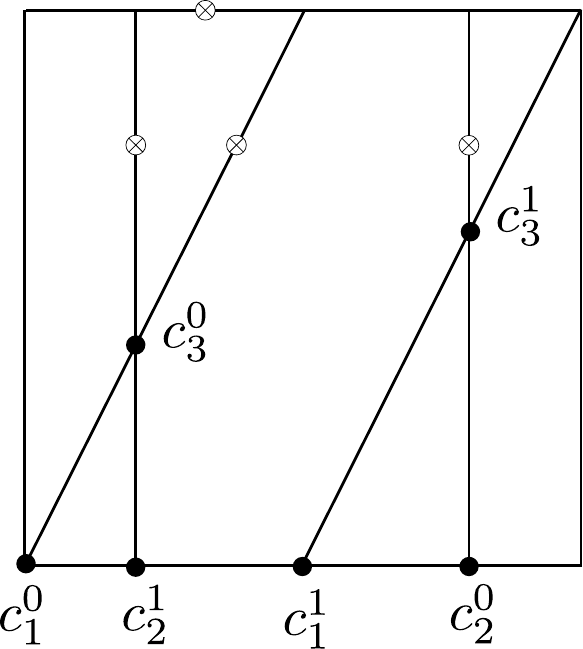}
    \caption{
    }
    \label{fig:Cones-L(1,2)-L(1,0)}
\end{figure}

\begin{step}
    $\mu^2(c_1,c_3^0) = 0$ is equivalent to $s([-q^{-x}M^{-1}]) = 0$,
    and $\mu^2(c_1,c_3^1) = 0$ is equivalent to $s([-q^xM]) = 0$.
\end{step}
\begin{proof}
    We will only show the statement for $c_3^0$, the other one being
    completely analogous. The signs with which the polygons
    encountered in the computation of $\mu^2(c_1,c_3^0)$ contribute
    will be determined according to the recipe described in
    \cite[Section 7]{Seidel--HMS-genus-2-surface} or \cite[Section
    2]{Lekili-Perutz--Fuk-torus-Dehn-surgery} (replacing $Spin$ by
    $Pin$). In particular, we think of the bounding $Pin$ structures
    on our Lagrangians as double covers which are trivialized except
    over one point where the two sheets are interchanged; these points
    are indicated by \raisebox{0.3mm}{$\scriptstyle{\otimes}$} in
    Figure \ref{fig:Cones-L(1,2)-L(1,0)}.
    
    The triangles contributing to $\mu^2(c_1,c_3^0)$ are precisely the
    images of the family of triangles $\Delta_n$, $n \in \mathbb Z$,
    in the universal cover $\mathbb R^2 \xrightarrow{\pi} T^2$
    described as follows: $\Delta_n$ has one vertex at $(x,2x)$, one
    at $(x,n)$, and one at $(\frac{1}{2}n,n)$, see Figure
    \ref{fig:triangles-univ-cov}. The first two points project to
    $c_3^0$ and $c_2^0$, respectively, while the third one projects to
    $c_1^0$ or to $c_1^1$, according to whether $n$ is even or
    odd. The area of $\Delta_n$ is
    \begin{equation*}
        A_n = \left(\frac{1}{2}(n-2x)\right)^2.
    \end{equation*}
    According to the recipe in \cite[Section
    7]{Seidel--HMS-genus-2-surface}, the sign with which
    $\pi(\Delta_n)$ contributes to $\mu^2(c_1,c_3^0)$ is $(-1)^{s_n +
        1}$ in the case at hand, where $s_n$ is the number of lifted
    \raisebox{0.3mm}{$\scriptstyle{\otimes}$} symbols encountered when
    travelling around the edges of $\Delta_n$. One sees easily that
    $s_n$ has the same parity as $n$, and therefore the sought-for
    sign is $(-1)^{n +1}$. 
    \begin{figure}[t]
        \centering
        \includegraphics[scale=0.6]{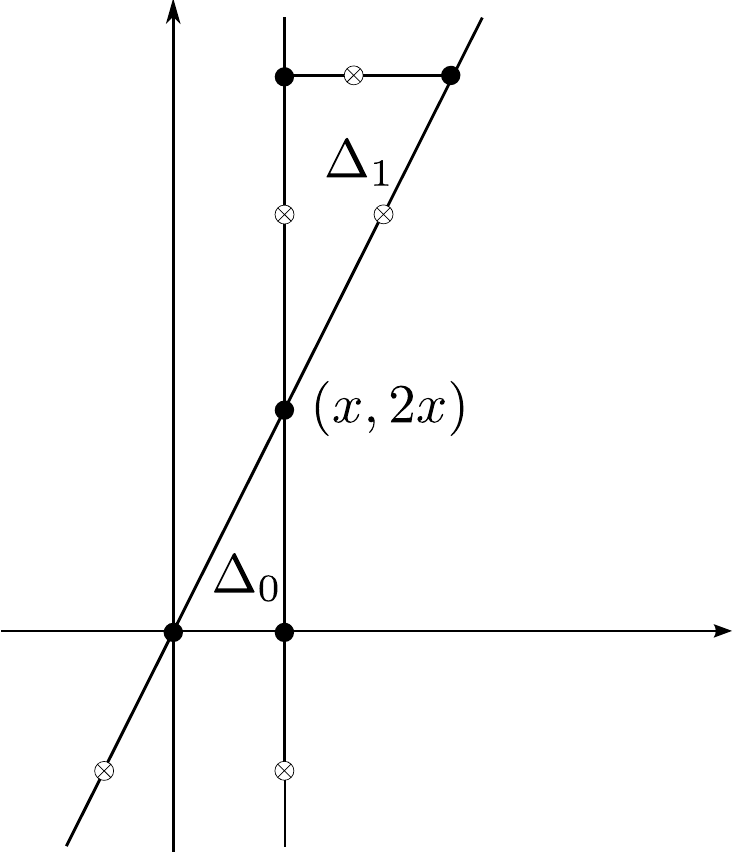}
        \caption{Triangles in the universal cover contributing to
            $\mu^2(c_1,c_3^0)$.}
        \label{fig:triangles-univ-cov}
    \end{figure}

    The remaining ingredient for the computation of $\mu^2(c_1,c_3^0)$
    is the parallel transport around the edges of
    $\pi(\Delta_n)$. Since the local systems on $Y_0$ and $Y_1$ are
    trivial, this is equal to the parallel transport in the local
    system $E_M^1$ along the vertical segment of (oriented) length
    $2x-n$ (the orientation of $L_{(0,-1)}$ points downwards); we
    denote this parallel transport map by $M^{2x-n}$.

    Assembling all these ingredients, we obtain
    \begin{equation*}
        \mu^2(c_1,c_3^0) = \left( \sigma^0 \sum_{n \in \mathbb Z}- M^{2x-2n}
            q^{(n-x)^2} + \sigma^1 \sum_{n \in \mathbb Z} M^{2x-(2n+1)}
            q^{(n + 1/2 - x)^2} \right) c_2^0.
    \end{equation*}
    Some basic arithmetic shows that the vanishing of this is
    equivalent to
    \begin{equation*}
        \sigma^0 \sum_{n \in \mathbb Z}(-q^{-x} M^{-1})^{2n} q^{n^2} +
        \sigma^1 \sum_{n \in \mathbb Z}(-q^{-x}M^{-1})^{2n+1}
        q^{(n+\frac{1}{2})^2} = 0,
    \end{equation*}
    which says nothing but $s([-q^{-x}M^{-1}]) = 0$.
\end{proof}

\begin{step}
    There exist a morphism $\widetilde c_1 \in
    \mathrm{hom}^0(Y_0,Y_1)$ such that $\mathrm{Cone}(\widetilde c_1)
    \cong Y_2$, where $Y_2$ is as defined in \eqref{eq:Y_2}.
\end{step}
\begin{proof}
    There exists a cobordism
    \begin{equation*}
        V: L_{(1,0)} \leadsto (L_{(0,-1),x},L_{(0,-1),-x},L_{(1,2)})
    \end{equation*}
    obtained from surgering $L_{(1,2)}$ with $L_{(0,-1),x}$ and
    $L_{(0,-1),-x}$, which results in $L_{(1,0)}$. It can be equipped
    with a local system restricting to $E_M^1$ and $E_{M^{-1}}^1$ on
    the two negative ends corresponding to $L_{(0,-1),x}$ and
    $L_{(0,-1),-x}$, and to trivial local systems on the other
    ends. Moreover, it can be equipped with a brane structure
    restricting to the standard brane structures of these Lagrangians
    as described in Section \ref{sec:comp-cobGr} (to verify this for
    the grading, use Lemma \ref{lem:gradings-cobsT2}). From this one
    can conlude the existence of a $\widetilde c_1$ with
    $\mathrm{Cone}(\widetilde c_1) \cong Y_2$ by ``bending'' the
    negative end of $V$ corresponding to $L_{(1,2)}$ such as to become
    positive, and then using a slight generalization of Theorem
    \ref{thm:cob-cone-decomp} for cobordisms with multiple ends on
    both sides.
\end{proof}

\begin{step}
    We have $\widetilde c_1 = c_1$ up to a non-zero factor.
\end{step}
\begin{proof}
    Since $\mathrm{Cone}(\widetilde c_1) \cong Y_2$, we infer that
    there is some non-zero morphism $c_3 = \eta^0c_3^0 + \eta^1 c_3^1
    \in \mathrm{hom}^1(Y_2,Y_0)$ such that $\mu^2(\widetilde c_1,c_3)
    = \eta^0 \mu^2(\widetilde c_1,c_3^0) + \eta^1\mu^2(\widetilde
    c_1,c_3^1) = 0$, again by \cite[Lemma
    3.7]{Seidel--Fukaya-Picard-Lefschetz-2008}.  This is equivalent to
    the individual vanishing of $\eta^0\mu^2(\widetilde c_1,c_3^0)$
    and $\eta^1\mu^2(\widetilde c_1,c_3^1)$, because the first is a
    multiple of $c_2^0$, while the second is a multiple of
    $c_2^1$. Since at least one of $\eta^0$ and $\eta^1$ is non-zero,
    we conclude that $\mu^2(\widetilde c_1, c_3^0) = 0$ or
    $\mu^2(\widetilde c_1,c_3^1) = 0$.
    
    Consider now the section $\widetilde s = \widetilde \sigma^0
    \theta^0 + \widetilde \sigma^1 \theta^1$ of $\mathcal O(2 P_0)$
    corresponding to $\widetilde c_1 = \widetilde \sigma^0 c_1^0 +
    \widetilde \sigma^1 c_1^1$. Assuming that $\mu^2(\widetilde
    c_1,c_3^0) = 0$, we conclude that $\widetilde s([-q^{-x}M^{-1}]) =
    0$ by the result of Step 1, and hence also $\widetilde s([-q^xM])
    = 0$, as $\widetilde s$ in a section of $\mathcal
    O(2P_0)$. (Assuming $\mu^2(\widetilde c_1,c_3^1) = 0$ would have
    led to the same conclusion.) 

    But this implies that $\widetilde s = s$ up to a non-zero factor,
    and hence $\widetilde c_1 = c_1$ up to a non-zero factor.
\end{proof}

Combing the results of Steps 2 and 3, we conclude that
\begin{equation*}
    \mathrm{Cone}(c_1) \cong Y_2,
\end{equation*}
as required. This finishes the prove of the lemma in the case that $s$
has two distinct zeros. 

We still have to argue that the cones of the four remaining morphisms,
for which the corresponding sections $s \in H^0(X;\mathcal O(2P_0))$
have a double zero, are as claimed. By a similar cobordism argument as
in the proof Step 2, we can infer that for every $x \in
\{0,\frac{1}{2}\}$ and every $M \in \{\pm 1\}$, there exist a morphism
$c_1 \in \mathrm{hom}^0(Y_0,Y_1)$ and an automorphism $a$ of
$(L_{(0,-1),x},E_M^1)$ such that $\mathrm{Cone}(c_1) \cong
\mathrm{Cone}(a)$. These morphisms $c_1$ must be the four remaining
ones, because we have already found the cones of all others to be
different from what we get here. Moreover, their cones must be
indecomposable, since the same is true on the algebraic side. But the
only indecomposable cone on an automorphism of $(L_{(0,-1),x},E_M^1)$
is $(L_{(0,-1),x},E_M^2)$.  \hfill$\square$

\subsection{Proof of Proposition \ref{prop:recov-mirr-funct}.}
Proposition \ref{prop:recov-mirr-funct} (i) for $h = 1$ has already
been proven. We break the rest of the proof up into several steps. The
common strategy is to exhibit for a given sheaf $\mathcal{F}$ a cone
decomposition in $D^b(X)$ with linearization consisting of objects for
which we already know Proposition \ref{prop:recov-mirr-funct} holds,
and with the property that $\mathcal F$ is \emph{determined} by the
cone decomposition.\footnote{Here and in the following, we should
    often add ``up to isomorphism'' to be really precise.} We then
construct, using iterated surgery, an object of $D^\pi \FukExt(T^2)$
which has a cone decomposition mirror to the one of $\mathcal{F}$ and
therefore must be mirror to $\mathcal{F}$.

\setcounter{step}{0}
\begin{step}
    Proposition \ref{prop:recov-mirr-funct} (ii) holds for all line
    bundles $\mathcal E \in \mathcal V(1,d)$, $d \in \mathbb Z$.
\end{step}
\begin{proof}
    We can write $\mathcal E = \mathcal O((d+1)P_0-Q)$ for some $Q \in
    X$, and there is an exact triangle
    \begin{equation*}
        \mathcal O_Q[-1] \to \mathcal E \to \mathcal O((d+1)P_0) \to \mathcal O_Q
    \end{equation*}
    in $D^b(X)$.
    Since $\mathrm{Hom}_{D^b(X)} (\mathcal O((d+1)P_0), \mathcal O_Q)$
    is one-dimensional, $\mathcal E$ is in fact the \emph{only}
    indecomposable object fitting into an exact triangle with the
    other two middle fixed. To find its mirror object, it suffices
    therefore to exhibit an exact triangle in $D^\pi \FukExt(T^2)$ in
    which two rightmost objects are mirror to $\mathcal O((d+1)P_0)$
    and $\mathcal O_Q$.

    We know already that these mirror objects are $L_{(1,-(d+1))}$ and
    $(L_{(0,-1),x},E_M^1)$ for certain $x \in \mathbb R/\mathbb Z$ and
    $M \in S^1\Lambda$. Now there exists some $x' \in \mathbb
    R/\mathbb Z$ such that surgering $L_{(0,-1),x}$ and
    $L_{(1,-d),x'}$ produces $L_{(1,-(d+1))}$, up to Hamiltonian
    isotopy. 
    After equipping the corresponding cobordism $V: L_{(1,-(d+1))}
    \leadsto (L_{(0,-1),x},L_{(1,-d),x'})$ with an appropriate brane
    structure and a local system
    , Theorem \ref{thm:cob-cone-decomp} yields an exact
    triangle
    \begin{equation*}
        (L_{(0,-1),x},E^1_{M})[-1] \to (L_{(1,-d),x'},E^1_{M^{-1}}) \to L_{(1,-(d+1))} \to  (L_{(0,-1),x},E^1_{M})
    \end{equation*}
    We conclude that $\Phi(\mathcal E) = (L_{(1,-d),x'},E^1_{M^{-1}})$
    up to isomorphism, as required.
    \end{proof}

    \begin{step}
        Proposition \ref{prop:recov-mirr-funct} holds for every
        $\mathcal E \in \mathcal V(r, d)$ whenever $\mathrm{gcd}(r,d)
        = 1$.
    \end{step}
    \begin{proof}
        Theorem \ref{thm:Atiyah-Thm3} implies that there exists an
        integer $n$ such that $\mathcal E$ fits into an exact triangle
        \begin{equation*}
            \mathcal L [-1] \to (\mathcal L')^{\oplus r-1} \to \mathcal E \to
            \mathcal L
        \end{equation*}
        with $\mathcal L = (\mathrm{det}\, \mathcal E)((r-1)n)$ and
        $\mathcal L' = \mathcal O_X(-n)$, where $-(k)$ denotes
        tensoring by the $k^{th}$ power of the hyperplane bundle. 
        Hence $\mathcal E$ admits an iterated cone decomposition with
        linearization $(\mathcal L, \mathcal L', \dots, \mathcal L')$,
        where $\mathcal L'$ appears $r-1$ times (cf. Lemma
        \ref{lem:iterated-icds}).

        We claim that $\mathcal E$ is actually the only indecomposable
        object of $D^b(X)$ admitting a cone decomposition with this
        linearization. Suppose that $\widetilde{\mathcal E}$ is
        another indecomposable object with the same property. It
        follows that $[\widetilde{ \mathcal E}] = [\mathcal E] \in
        K_0(X)$ by Lemma \ref{lem:linearization-and-K0}, and hence
        $\mathrm{det}\,\widetilde{\mathcal E} = \mathrm{det}\,\mathcal
        E$ and $\mathrm{rk}\,\widetilde{\mathcal E} =
        \mathrm{rk}\,\mathcal E$, because $(\mathrm{det},
        \mathrm{rk}): K_0(X) \to \mathrm{Pic}(X) \oplus \mathbb Z$ is
        an isomorphism (see \cite{Hartshorne--AG}). Since
        $\widetilde{\mathcal E}$ is indecomposable, we conclude that
        $\widetilde{\mathcal E} \in \mathcal V(r,d)$. By Theorem
        \ref{thm:Atiyah-Thm6}, the condition $\mathrm{gcd}(r,d)= 1$
        implies that $\mathrm{det}: \mathcal V(r,d) \to \mathcal
        V(1,d) \subset \mathrm{Pic}(X)$ is bijective, and hence
        $\widetilde{\mathcal E} = \mathcal E$.

        To find $\Phi(\mathcal E)$, it is therefore sufficient to
        construct an object of $D^\pi \FukExt(T^2)$ admitting a cone
        decomposition with linearization mirror to $(\mathcal L,
        \mathcal L',\dots,\mathcal L')$. The line bundles there have
        degrees $\mathrm{deg}\,\mathcal L = d + 3n(r-1)$ and
        $\mathrm{deg}\,\mathcal L' = -3n$, since the hyperplane bundle
        has degree $3$. 
        We conclude, using the previous step, that their mirror
        objects are of the form 
        \begin{equation*}
            \Phi(\mathcal L) = (L,E) \quad \text{and} \quad \Phi(\mathcal L') = (L',E'),
        \end{equation*}
        with Lagrangians $L$ and $L'$ of slopes $(1,-d-3n(r-1))$ and
        $(1,3n)$ that are equipped with rank one local systems $E$
        and $E'$. Starting with $L$ and iteratively surgering $r-1$
        times with $L'$ leads to a sequence of Lagrangians
        $N_0=L,N_1,\dots,N_{r-1}$ of slopes
        \begin{equation*} 
            [N_j] = (1+j,-d-3n(r-1-j)) \in H_1(T^2;\mathbb Z), \quad j = 0,\dots,r-1.
        \end{equation*}
        (Note that $N_j$ might have multiple components, and a priori
        these might have slopes equal to that of $L'$, which would be
        problematic. 
        But the only way this can happen is that $N_j$ has $1+j$
        components of slope $(1,(-d-3n(r-1-j))/(1+j))$; to exlude this
        being equal to $(1,3n) = [L']$, we assume that $n$ has been
        chosen such that $-d < 3nr$. This is no restriction, since any
        sufficiently big $n$ leads to an exact triangle as above.)

        To this sequence of surgeries corresponds a sequence of
        cobordisms 
        \begin{equation*}
            V_j: N_j \leadsto (L, \underbrace{L',\dots,L'}_{j
                ~\text{times}}), \quad j = 0,\dots,r-1.
        \end{equation*}
        These can be equipped with brane structures and rank one local
        systems that extend those on the negative ends and determine a
        brane structure and a rank one local system $E_j$ on the
        positive end.  Theorem \ref{thm:cob-cone-decomp} applied to
        $V_{r-1}$ says that there is an iterated cone decomposition of
        $(N_{r-1},E_{r-1})$ with linearization
        \begin{equation*}
            ((L,E),(L',E'),\dots,(L',E')).
        \end{equation*}
        This is mirror to the one for $\mathcal E$, and we infer
        $\Phi(\mathcal E) \cong (N_{r-1},E_{r-1})$. That proves the
        claimed statement, since $N_{r-1}$ has slope $(r,-d)$.
    \end{proof}

\begin{step}
    \label{step:gcd-greater-one} Proposition
    \ref{prop:recov-mirr-funct} holds for all $\mathcal E \in \mathcal
    V(hr,hd)$ with $\mathrm{gcd}(r,d) = 1$ and $h \geq 1$, and for all
    indecomposable skyscraper sheaves with stalks of rank $h \geq 1$.
\end{step}
\begin{proof}
    By Lemma \ref{lem:Atiyah-lem26}, every indecomposable vector
    bundle $\mathcal E$ on $X$ is of the form $E_{\mathcal L_1}(hr,hd)
    \otimes \mathcal L_0$ for some $(r,d)$ with $\mathrm{gcd}(r,d) =
    1$, $h \in \mathbb N$, and certain line bundles $\mathcal
    L_0,\mathcal L_1$ of degrees zero and one (see Appendix
    \ref{app:vect-bundl-ellipt} for the notation). Moreover, for every
    $h \geq 1$ and every $Q \in X$, there is a unique indecomposable
    skyscraper sheaf $\mathcal O_{hQ}$ with stalk of rank $h$
    supported at $Q$, and every indecomposable skyscraper sheaf is of
    this form.

    Denote by $\mathcal{Y}_h$ either the vector bundle $E_{\mathcal
        L_1}(hr,hd) \otimes \mathcal L_0$ for fixed $r,d,\mathcal L_0$
    and $\mathcal L_1$, or the skyscraper sheaf $\mathcal O_{hQ}$ for
    fixed $Q$. We will prove that for any such choice, the mirror
    objects of the $\mathcal{Y}_h$, $h \geq 1$, are of the form
    $(L,E^h)$, with $L$ a Lagrangian of slope $(r,d)$, or $(0,-1)$,
    and $E^h$ an indecomposable local system of rank $h$ over $L$
    such that there exists a short exact sequence
    \begin{equation*}
        0 \to E^1 \to  E^{h+1} \to E^h \to 0
    \end{equation*}
    of local systems on $L$ (which in this case is equivalent to
    saying that the unique eigenvalue of the monodromy is the same for
    all of them). 
    The proof will be by induction on $h$. The claim for $h = 1$ is
    what was proven in the previous steps.

    In both of the above cases and for every $h \in \mathbb N$, there
    is an exact triangle
    \begin{equation*}
        \mathcal{Y}_h[-1] \to \mathcal{Y}_1 \to \mathcal{Y}_{h+1} \to
        \mathcal Y_h
    \end{equation*}
    in $D^b(X)$. For $\mathcal{Y}_h = E_{\mathcal L_1}(hr,hd) \otimes
    \mathcal L_0$, this comes from the short exact sequence obtained
    by tensoring the short exact sequence in Theorem
    \ref{thm:Atiyah-Thm5}(i) with $E_{\mathcal L_1}(hr,hd) \otimes
    \mathcal L_0$ (also cf. \ref{lem:Atiyah-lem24}). For
    $\mathcal{Y}_h = \mathcal O_{hQ}$, it comes from the short exact
    sequence $0 \to \mathcal O_Q \to \mathcal O_{(h+1)Q} \to \mathcal
    O_{hQ} \to 0$.  Moreover, we have
    \begin{equation*}
        \mathrm{dim}~ \mathrm{Ext}^1 (\mathcal{Y}_h,\mathcal{Y}_1) = 1
    \end{equation*}
    in both cases (to see this in the first case, use Serre duality
    and \cite[Lemma 22]{Atiyah--VecBunEllCurve--1957}.) Hence
    $\mathcal{Y}_{h+1}$ is the only indecomposable object of $D^b(X)$
    that can arise as a cone on a morphism $\mathcal{Y}_h[-1] \to
    \mathcal{Y}_1$.

    By inductive assumption, the mirror images of $\mathcal Y_1$ and
    $\mathcal Y_h$ share the same Lagrangian brane that's equipped
    with indecomposable local systems $E^1, E^h$ of ranks $1$ and $h$,
    and whose monodromies have the same eigenvalue (because there
    exists a morphism of local systems $E^1 \to E^h$). 
    There's a unique indecomposable local system $E^{h+1}$ of rank
    $h+1$ fitting into a short exact sequence $0 \to E^{1} \to E^{h+1}
    \to E^{h}\to 0$ of local systems on $L$ (namely, the unique
    indecomposable local system whose monodromy has this
    eigenvalue). Hence there exists an exact triangle
    \begin{equation*}
        (L, E^h)[-1] \to (L,E^1) \to  (L,E^{h+1}) \to (L,E^h)
    \end{equation*}
    in $D \FukExt(T^2)$ by Proposition
    \ref{prop:ses-ls-exact-triangle}. We conclude that $\Phi(\mathcal
    Y_{h+1}) = (L,E^{h+1})$ by the same argument as in the previous
    steps. This is as required.
\end{proof}

This ends the proof of Proposition
\ref{prop:recov-mirr-funct}. \hfill$\square$

\section{Proofs of the main theorems}
\label{sec:isomorphism}

We will first prove Theorem \ref{thm:MainThmExtended} and then deduce
Theorems \ref{thm:MainThm} and \ref{thm:MainThm-SES}.

\subsection{Proof of Theorem \ref{thm:MainThmExtended}.}
\label{sec:proof-theor-refthm:m}
The cobordism group $\cobGrExt(T^2)$ appearing in Theorem
\ref{thm:MainThmExtended} is defined as described in Section
\ref{sec:relations-from-local-systems} with $\LagsExt (T^2) =
\mathrm{Ob}\,\FukExt(T^2)$, i.e., it has as generators Lagrangian
branes with local systems and relations coming from cobordisms with
vanishing Maslov class carrying compatible local systems and gradings,
as well as additional relations induced by short exact sequences of
local systems.

Let $X$ be the Tate curve mirror to $T^2$, and denote by
\begin{equation*}
    K_0(X) = K_0(D^b(X))
\end{equation*}
the Grothendieck group of its derived category of coherent
sheaves. The mirror functor $\Phi: D^b(X) \to D^\pi \FukExt (T^2)$ of
Theorem \ref{thm:HMS-for-T2} induces an isomorphism of Grothendieck
groups as it is an equivalence of triangulated categories. Since
moreover the inclusion $D\FukExt(T^2) \hookrightarrow
D^\pi\FukExt(T^2)$ is an equivalence by Corollary
\ref{cor:DF-is-equiv-to-DpiF}, and hence also induces an isomorphism
of Grothendieck groups, we obtain an isomorphism
\begin{equation*}
    K_0(D \FukExt(T^2)) \xrightarrow{\cong} K_0(X).
\end{equation*}

Denote by $\underline{\mathrm{Coh}}\,X \subset D^bX$ the full
subcategory consisting of direct sums of shifted sheaves; as mentioned
before, $\underline{\mathrm{Coh}}\,X$ is equivalent to
$D^bX$. Consider the map
\begin{equation*}
    F_\Phi: \langle  D^b(X)\rangle \to \langle \LagsExt(T^2) \rangle
\end{equation*}
induced by $\Phi$, in the sense that it first replaces an arbitrary
object of $D^bX$ by an isomorphic object of
$\underline{\mathrm{Coh}}\,X$, and then applies the mirror functor to
get a sum of Lagrangian branes with local systems. We will prove that
$F_\Phi$ descends to a well-defined group homomorphism $K_0(X) \to
\cobGrExt(T^2)$, thus completing the diagram
\begin{equation}
    \label{eq:cobGr-and-K_0}
    \begin{aligned}
        \xymatrix{ \cobGrExt(T^2) \ar[rr]^{\Theta^\sesLS\quad}&& K_0(D
            \FukExt(T^2))\ar[ld]^{\cong}\\
            &~~K_0(X)\ar@{-->}[lu]&
        }
    \end{aligned}
\end{equation}
where $\Theta^\sesLS: \cobGrExt(T^2) \to K_0(D \FukExt(T^2))$ is the
canonical surjective homomorphism from Proposition
\ref{prop:Cob-K-group-hom}. Once well-definedness is proven, it is
clear by construction that the composition of the two lower
homomorphisms in \eqref{eq:cobGr-and-K_0} provides a left-inverse for
$\Theta^\sesLS$, showing in particular that $\Theta^\sesLS$ is
injective. Given that we already have surjectivity, this will conclude
the proof of \mbox{Theorem
    \ref{thm:MainThmExtended}}. 

\subsubsection{Well-definedness.}
\label{sec:well-defin-gamma}
Denote by $R$ the set of $K_0$-relations among objects of $\langle
\underline{\mathrm{Coh}}\,X \rangle$, i.e., the kernel of $\langle
\underline{\mathrm{Coh}}\,X \rangle \to K_0(X)$. It follows from
Proposition \ref{prop:K-Abelian-K-derived} that the inclusion
$\underline{\mathrm{Coh}}\,X \hookrightarrow D^bX$ induces an
isomorphism $\langle \underline{\mathrm{Coh}}\,X \rangle/R \cong
K_0(X)$. To prove that $F_\Phi$ induces a map $K_0(X) \to
\cobGrExt(T^2)$, we must therefore show that $F_\Phi$ takes $R$ to
$R^\sesLS$, the kernel of $\langle \LagsExt(T^2) \rangle \to
\cobGrExt(T^2)$. The essential part is to prove that this is true for
the subset $ R_0 = R \cap \langle \mathrm{Ind\,Coh}\,X \rangle$, the
set of $K_0$-relations among indecomposable objects of
$\mathrm{Coh}\,X$.

\begin{prop}
    \label{prop:generators-of-relations}
    $R_0$ is generated by relations coming from short
    exact sequences of the following types:
    \begin{itemize}
    \item $0 \to \mathcal F \to \mathcal G \to 0$ with indecomposable
        $\mathcal F, \mathcal G \in \mathrm{Coh}\,X$;
    \item $0 \to \mathcal O(D-Q) \to \mathcal O(D) \to \mathcal O_Q
        \to 0$ with $D$ a divisor and $Q \in X$;
    \item $0 \to (\mathcal O_X^{\oplus r-1})(-n) \to \mathcal{E} \to
        (\mathrm{det} \,\mathcal{E})((r-1)n) \to 0$ with $\mathcal{E}
        \in \mathcal{V}(r,d)$ such that $\mathrm{gcd}(r,d) = 1$, and $n
        \in \mathbb Z$ such that $\mathcal{E}(n)$ is generated by
        global sections;
    \item $0 \to \mathcal{Y}_1 \to \mathcal{Y}_{h+1} \to
        \mathcal{Y}_h \to 0$ for $h \geq 1$.
    \end{itemize}
    Here $\mathcal{Y}_h$ denotes either a vector bundle
    of the form $E_{\mathcal L_0}(rh,dh) \otimes \mathcal L_1$ for
    some fixed $r,d$ with $\mathrm{gcd}(r,d) = 1$ and fixed line
    bundles $\mathcal L_0, \mathcal L_1$ of degree 0 resp. 1, or an
    indecomposable skyscraper sheaf $\mathscr O_{hq}$ for some $q \in
    X$. (The SESs of the second type are the obvious ones; those of
    the third type are as in Theorem \ref{thm:Atiyah-Thm3}; and those
    of the fourth type are as mentioned in Step 3 of the proof of
    Proposition \ref{prop:recov-mirr-funct}.)
\end{prop} 

\begin{proof}
    Denote by $S \subset \langle \mathrm{Ind\,Coh}\,X \rangle$ the
    subgroup generated by the relations induced by the short exact
    sequences of the types stated, and by $\llbracket \mathcal{F}
    \rrbracket$ the class of $\mathcal{F} \in \mathrm{Ind \,Coh}\,X$
    in $\langle \mathrm{Ind\,Coh}\,X \rangle / S$. We will show that
    \begin{equation}
        \label{eq:det-rk-equality}
        \llbracket \mathcal{F} \rrbracket = \llbracket
        \mathrm{det}\,\mathcal{F} \rrbracket + (\mathrm{rk}\,\mathcal{F}-1)\llbracket
        \mathcal{O}_X \rrbracket.
    \end{equation}
    This implies immediately that the canonical map
    $\langle\mathrm{Ind\,Coh}\,X \rangle / S \to K_0(X)$, which is
    easily seen to be surjective using Proposition
    \ref{prop:K-Abelian-K-derived}, is also injective: Namely, since
    $(\mathrm{det},\mathrm{rk}): K_0(X) \to \mathrm{Pic}(X) \oplus
    \mathbb Z$ is an isomorphism, the equality $[\mathcal{F}] =
    [\mathcal{G}]$ in $K_0(X)$ is equivalent to
    $\mathrm{det}\,\mathcal{F} = \mathrm{det}\,\mathcal{G}$ and
    $\mathrm{rk}\,\mathcal{F} = \mathrm{rk}\,\mathcal{G}$, from which
    $\llbracket \mathcal{F} \rrbracket = \llbracket \mathcal{G}
    \rrbracket$ follows by \eqref{eq:det-rk-equality}. This implies
    that the described relations generate all of $R_0$.
    
    We first claim that every sheaf can be written as a sum or
    difference of line bundles in $\langle \mathrm{Ind\,Coh}\,X
    \rangle / S$. Now an indecomposable sheaf on $X$ is either a
    vector bundle of the form $E_{\mathcal L_0}(rh,dh) \otimes
    \mathcal L_1$ for certain $r,d$ with $\mathrm{gcd}(r,d) = 1$ and
    $h \geq 1$, or a skyscraper sheaf $\mathscr O_{hq}$ for some $q
    \in X$ and $h \geq 1$. Using the relations of the last type, we
    can inductively reduce to considering the case $h = 1$, i.e., the
    case of vector bundles $\mathcal{E} \in \mathcal{V}(r,d)$ with
    $\mathrm{gcd}(r,d) = 1$ or of skyscraper sheaves $\mathcal{O}_Q$.
    Now the relations of second and third types show that these
    satisfy
    \begin{equation*}
        \llbracket \mathcal{E} \rrbracket = (r-1)\llbracket \mathcal
        O_X(-n) \rrbracket + \llbracket
        \mathrm{det}\,\mathcal{E}((r-1)n) \rrbracket,\quad \llbracket
        \mathcal O_Q \rrbracket = \llbracket \mathcal{O}(D) \rrbracket -
        \llbracket \mathcal{O}(D-Q) \rrbracket,
    \end{equation*}
    and since the classes on the right-hand sides are those of lines
    bundles, our claim is shown.

    We now prove that
    \begin{equation}
        \label{eq:relations-for-line-bundles}
        \llbracket \mathcal{L} \rrbracket + \llbracket \mathcal{L}'
        \rrbracket = \llbracket \mathcal{L} \otimes \mathcal{L}'
        \rrbracket + \llbracket \mathcal{O}_X \rrbracket, \quad \llbracket \mathcal{L} \rrbracket - \llbracket \mathcal{L}'
        \rrbracket = \llbracket \mathcal{L} \otimes \mathcal{L}'^{-1}
        \rrbracket
    \end{equation}
    for all lines bundles $\mathcal{L}$ and $\mathcal{L}'$. Observe
    that every line bundle is of the form $\mathcal O(D)$ for a
    divisor $D = \sum_{i=1}^n P_i - \sum_{j=1}^m Q_j$. Using the
    relations of the second type in the statement of the lemma, one
    obtains inductively that
    \begin{equation*}
        \llbracket \mathcal O(D) \rrbracket = \llbracket \mathcal{O}_X
        \rrbracket + \sum_{i=1}^n \llbracket \mathcal{O}_{P_i}
        \rrbracket - \sum_{j=1}^m \llbracket \mathcal{O}_{Q_j} \rrbracket.
    \end{equation*}
    The identities in \eqref{eq:relations-for-line-bundles} follow
    from this because $\mathcal{O}(D) \otimes \mathcal{O}(D') =
    \mathcal{O}(D + D')$ and $\mathcal{O}(D) \otimes
    \mathcal{O}(D')^{-1} = \mathcal{O}(D - D')$. 

    Equation \eqref{eq:det-rk-equality} now follows easily for all
    $\mathcal{F} \in \mathrm{Ind \,Coh}\,X$. As we have seen, we can
    write $\llbracket \mathcal{F} \rrbracket = \sum_{i=1}^k s_i
    \llbracket \mathcal{L}_i \rrbracket$ with lines bundles $\mathcal
    L_i$ and $s_i \in \{\pm 1\}$, and then $\mathrm{det}\,\mathcal{F}
    = \bigotimes_{i=1}^k \mathcal{L}_i^{s_i}$ and
    $\mathrm{rk}\,\mathcal{F} = \sum_{i=1}^k s_i$. Applying
    \eqref{eq:relations-for-line-bundles} inductively yields
    \eqref{eq:det-rk-equality}.
\end{proof}

Generators of $R_0$ of the first type listed in Proposition
\ref{prop:generators-of-relations} identify isomorphic indecomposable
sheaves; their mirror images are Hamiltonian isotopic curves with
isomorphic local systems, and we know that there exists a cobordism
between these. Cobordisms corresponding to the generators of $R_0$ of
the other types were constructed in the proof of Proposition
\ref{prop:recov-mirr-funct}. This implies that $F_\Phi$ takes $R_0$ to
$R^\sesLS$.

The remaining relations in $R$ come from taking direct sums and
shifting; that is, we can write $R = R_0 + R_1$ where $R_1$ is
generated by elements of the form $\mathcal F \oplus \mathcal G -
(\mathcal F + \mathcal G)$ and $\mathcal F[1] + \mathcal F$. That
$F_\Phi$ takes these generators and hence all of $R_1$ to $R^\sesLS$
is clear by definition for those of type 'direct sum', and follows for
those of type 'shift' from the fact that $(L[1],L)$ is null-cobordant
for every Lagrangian brane $L$: A null-cobordism is given by $V =
\gamma \times L$ equipped with a suitable grading, where $\gamma \subset
\mathbb R^2$ is a curve with two negative ends.

This finishes the proof of Theorem \ref{thm:MainThmExtended}.
\hfill$\square$

\subsection{Proof of Theorem \ref{thm:MainThm}.}
\label{sec:proof-theorem-b}
Consider the diagram
\begin{equation*}
    \begin{aligned}
        \xymatrix@-8pt{
            \cobGr(T^2) \ar[r] \ar[d]_{\Theta} & \cobGrExt(T^2)
            \ar[d]^{\Theta^\sesLS}\\
            K_0(D \Fuk(T^2)) \ar[r] & K_0(D \FukExt(T^2))
            }
    \end{aligned}
\end{equation*}
in which the upper arrow is the canonical map induced by $\Lags(T^2)
\hookrightarrow \,\LagsExt(T^2)$, and where the lower arrow is the map
induced by the inclusion $D\Fuk(T^2) \hookrightarrow D
\FukExt(T^2)$. The square commutes because both compositions take the
class of any given Lagrangian $L$ in $\cobGr(T^2)$ to the class of $L$
in $K_0(D\FukExt(T^2))$. Note that the upper horizontal map is
injective, as it is a section of the group homomorphism
$\cobGrExt(T^2) \to \cobGr(T^2)$ induced by $(L,E) \mapsto
\mathrm{rk}(E) L$ (that this is well-defined is immediate from the
definition of the relations in both groups, see Section
\ref{sec:lagr-cobord-group}). Since $\Theta^\sesLS$ is already known to be
an isomorphism, this implies that $\Theta$ is an isomorphism and hence
concludes the proof of Theorem \ref{thm:MainThm}. \hfill $\square$

\subsection{Proof of Theorem \ref{thm:MainThm-SES}}
\label{sec:proof-theorem-SES}
In view of Proposition \ref{prop:an-exact-sequence}, we are left with
showing that the map $\zeta: \mathbb R/\mathbb Z \to \cobGr(T^2)$
defined in Section \ref{sec:cylinders} is injective. To see this,
consider the composition
\begin{equation*}
    \mathbb R/\mathbb Z \xrightarrow[]{\zeta}
    \cobGr(T^2) \xrightarrow[]{\Theta} K_0(D\Fuk(T^2)) \to
    K_0(D\FukExt(T^2)) \to K_0(X).
\end{equation*}
Recall that $\zeta(x) = [L_{(0,1),x}] - [L_{(0,1)}] = [L_{(0,-1)}] -
[L_{(0,-1),x}]$. By what we know about the action of the mirror
functor (see the proof of Proposition \ref{prop:recov-mirr-funct}),
the composition hence takes $x \mapsto [\mathcal O_{[-q^0]}] -
[\mathcal O_{[-q^{x}]}] \in K_0(X)$, which is zero if and only if $x =
0 \in \mathbb R/\mathbb Z$. So the entire composition is injective,
and therefore $\zeta$ is. \hfill $\square$

\appendix
\section{Iterated cone decompositions and $K_0$}
\label{sec:a_infty-triang-categ}
\subsection{Triangulated categories.}
\label{sec:triang-categ}

A \emph{triangulated category} $\mathscr D$ is an additive category
equipped with an additive autoequivalence $S: \mathscr D \to \mathscr
D$ called the \emph{shift functor}, and a set of \emph{exact
    triangles} $ X \xrightarrow{f} Y \xrightarrow{g} Z \xrightarrow{h}
S(X)$. These data are required to satisfy a list of axioms, for which
we refer to \cite{Weibel--Intro-to-Hom-Alg}. The most relevant for us
is that every morphism $f: X \to Y$ in $\mathscr D$ can be
completed to such an exact triangle. The object $Z$ is then determined
up to isomorphism, and we call it a \emph{cone} on the morphism
$f$. Moreover, we write $X[n] = S^nX$ and $f[n] = S^n(f)$ for the
effect of iterates of the shift functor $S$ on objects and
morphisms. This is in reminiscence of the homotopy category of
complexes $K(\mathscr A)$ over an Abelian category $\mathscr A$, which
is the prototypical example of a triangulated category.

\subsection{Generation and iterated cone decompositions.}
\label{sec:gener-iter-cone}
Given a full subcategory of a triangulated category $\mathscr D$ with
objects a collection $\{X_i ~|~ i \in I\}$, one can consider the full
subcategory consisting of all objects that are cones on morphisms
between the $X_i$. Iterating this construction, i.e.\ including in
each step all cones on morphisms between the previously constructed
objects, one ends up with the subcategory of $\mathscr D$
\emph{generated} by the $X_i$. 

In the other direction, one can ask if and how an object $X$ of
$\mathscr D$ can be constructed as an iterated cone on morphisms
between other objects. The following notion is useful to formalize
this.

\begin{defn} Let $\mathscr D$ be a triangulated category and let $X
    \in \mathscr D$. An \emph{iterated cone decomposition of $X$} is a
    sequence of exact triangles
    \begin{equation*}
        C_{i-1}[-1] \to X_i \to C_i \to C_{i-1}, \quad i = 1, \dots, k,
    \end{equation*}
    with objects $C_0,C_1,\dots, C_k \in \mathscr D$ that satisfy $C_0
    = 0$ and $C_k = X$. The tuple $(X_1,\dots,X_k)$ is called the
    \emph{linearization} of the cone decomposition.
\end{defn}

\begin{rmk}
    This definition is adapted to the cohomological conventions we use
    and therefore differs from the one in \cite[Section
    2.6]{Biran-Cornea--Lag-Cob-II}, where homological conventions are
    used.
\end{rmk}

Iterated cone decompositions can themselves be iterated and are
well-behaved with respect to that in the sense of the following lemma
(cf.\ the composition in Biran--Cornea's category $T^S\mathscr D$ of
iterated cone decompositions \cite[Section
2.6]{Biran-Cornea--Lag-Cob-II}). 

\begin{lem}
    \label{lem:iterated-icds}
    Suppose that $X \in \mathscr D$ admits an iterated
    cone decomposition with linearization $(X_1,\dots,X_k)$, and that
    one of the $X_h$, $1 \leq h \leq k$, admits an iterated cone
    decomposition with linearization $(X_h^1,\dots,X_h^\ell)$. Then
    $X$ admits an iterated cone decomposition with linearization
    \begin{equation*}
        (X_1,\dots,X_{h-1},X_h^1,\dots,X_h^\ell,X_{h+1},\dots,X_k).
    \end{equation*}
\end{lem}

\subsection{Grothendieck groups.}
\label{sec:k-groups}
Let $\mathscr D$ be a triangulated category. The Grothen\-dieck group
$K_0(\mathscr D)$ is defined as the quotient $ K_0(\mathscr D) =
\langle \mathrm{Ob}\,\mathscr D \rangle / R$ of the free Abelian group
generated by the objects of $\mathscr D$ by the subgroup $R$ generated
by all expressions $X - Y + Z \in \langle \mathrm{Ob}\,\mathscr D
\rangle$ such that there exists an exact triangle $X \to Y \to Z \to
X[1]$. The following lemma is straightforward.
\begin{lem}
    \label{lem:linearization-and-K0}
    Suppose that $X \in \mathscr D$ admits an iterated cone
    decomposition with linearization $(X_1,\dots,X_k)$. Then $[X] =
    [X_1] + \dots + [X_k]$ in $K_0(\mathscr D)$.
\end{lem}

One can also define the Grothendieck group $K_0(\mathscr A)$ of an
Abelian category $\mathscr A$, by starting with the free Abelian group
on objects and imposing a relation $[A] + [C] = [B]$ for every short
exact sequence $0 \to A \to B \to C \to 0$. Recall the canonical
inclusion $\mathscr A \hookrightarrow D^b(\mathscr A)$, which on
objects is given by viewing $X \in \mathrm{Ob}\,\mathscr A$ as a
complex concentrated in degree zero. The following statement is
well-known.

\begin{prop}
    \label{prop:K-Abelian-K-derived}
    The canonical inclusion $\mathscr A \hookrightarrow D^b(\mathscr
    A)$ induces an isomorphism $K_0(\mathscr A) \cong K_0(D^b\mathscr
    A)$.
\end{prop}

\section{Exact triangles from SESs of local systems}
\label{app:SESs-LSs}
Consider a symplectic manifold $(M,\omega)$ for which we can define
the Fukaya category $\FukExt(M)$ with gradings and signs as outlined
in Section \ref{sec:Fuk-T2} such that the objects are Lagrangian
branes with local systems. There appears to be no reference in the
literature for the following statement.

\begin{prop}
    \label{prop:ses-ls-exact-triangle}
    Let $L^\deco$ be a Lagrangian brane 
    and let $0 \to E' \to E \to E'' \to 0$ be a short
    exact sequence of local systems on $L$. Then there exists an exact
    triangle
    \begin{equation*}
        (L^\deco,E'')[-1] \to (L^\deco,E') \to (L^\deco,E) \to (L^\deco,E'')
    \end{equation*}
    in $D \FukExt(M)$.
\end{prop}

We view local systems as assignments of vector spaces and parallel
transport maps as described in Section \ref{sec:local-systems}. By a
short exact sequence of local systems $0 \to E' \xrightarrow{i} E
\xrightarrow{p} E'' \to 0$ we mean a family of short exact sequences
of vector spaces
\begin{equation*}
    0 \to E'_x \xrightarrow{i_x} E_x \xrightarrow{p_x}
    E''_x \to 0
\end{equation*}
for every $x \in L$ such that the $i_x$ and $p_x$ define morphisms of
local systems, i.e., such that they commute with parallel transport
maps. For the following proof, we choose a splitting
\begin{equation*}
    0~ \to ~E'_x~ 
    \raisebox{-2pt}{$\stackrel{\xrightarrow{i_x}}{\xleftarrow[q_x]{}}
        $}
    ~E_x~
    \raisebox{-2pt}{$\stackrel{\xrightarrow{p_x}}{\xleftarrow[j_x]{}}
        $}
    ~E''_x~ \to ~0
\end{equation*}
for every $x \in L$, that is, maps $q_x: E_x \to E'_x$ and $j_x: E''_x
\to E_x$ such that $q_x \circ i_x = \mathrm{id}_{E'_x}$, $p_x \circ
j_x = \mathrm{id}_{E''_x}$ and $i_x\circ q_x + j_x \circ p_x =
\mathrm{id}_{E_x}$. Note that the $j_x$ and $q_x$ do generally not
define morphisms of local systems, unless the short exact sequence
splits globally.

In the proof of Proposition \ref{prop:ses-ls-exact-triangle}, it will
be convenient to model the relevant morphism spaces in $\FukExt(M)$ as
spaces of Morse cochains with coefficients in local systems. This is
possible under certain conditions on $M$ and $L$, the fundamental
example being that $M$ is a cotangent bundle and $L$ is an exact
Lagrangian (from which the case of interest in this
paper follows immediately by a covering
argument).

We will adapt a construction used in
\cite{Abouzaid--A-top-model-Fuk-plumbings} (which goes back to
\cite{Fukaya-Oh--Zero-loop-open-strings}) and consider an
$A_\infty$-category $\mathscr M(L)$ defined as follows: The objects of
$\mathscr M(L)$ are all those objects of $\FukExt(M)$ whose underlying
Lagrangian is $L$. For the morphism spaces, we fix a Morse function
$f: L \to \mathbb R$ and define
\begin{equation*}
    \mathrm{hom}_{\mathscr M(L)}^i (E_0,E_1) = \bigoplus_{\substack{x
            \in \mathrm{Crit}\,f \\\vert x
            \vert = i}} \mathrm{Hom}(E_{0,x},E_{1,x}),
\end{equation*}
where $\vert \cdot \vert$ denotes the Morse index; here we denote
objects of $\mathscr M(L)$ simply by their local system. The
$A_\infty$-operations $\mu^d_{\mathscr M(L)}$ are defined by
considering rigid perturbed gradient flow trees with vertices at
critical points, and summing up parallel transport maps in the
relevant local systems along edges of these trees. See
\cite{Abouzaid--A-top-model-Fuk-plumbings} for the description of the
relevant moduli spaces.

By adapting the arguments in
\cite{Abouzaid--A-top-model-Fuk-plumbings}, one can show (under
certain conditions, as indicated above) that there is an
$A_\infty$-quasi-isomorphism
\begin{equation*}
    \FukExt(L) \to \mathscr M(L),
\end{equation*}
where $\FukExt(L)$ denotes the full $A_\infty$-subcategory of
$\FukExt(M)$ consisting of objects whose underlying Lagrangian is
$L$. It will therefore suffice to prove that there exists an exact
triangle of the claimed form in $H^0(Tw \mathscr M(L))$, where $Tw
\mathscr M(L)$ denotes the category of twisted complexes over
$\mathscr M(L)$ which we use to model the triangulated closure of that
category (see \cite[Section
(3l)]{Seidel--Fukaya-Picard-Lefschetz-2008}).

\begin{proof}[Proof of Proposition \ref{prop:ses-ls-exact-triangle}]
    We assume that the Morse function $f: L \to \mathbb R$ defining
    the morphism spaces in $\mathscr M(L)$ has a single local minimum
    $x_0 \in L$.  As for notation, we write $\pi_\gamma',\pi_\gamma,
    \pi_\gamma''$ for the parallel transport in $E',E,E''$ along a
    path $\gamma$, and we denote by $\overline \gamma$ the path
    obtained by reversing $\gamma$.

    Let $c_2 \in \mathrm{hom}^0(E',E) =
    \mathrm{Hom}(E_{x_0}',E_{x_0})$ and $c_3 \in \mathrm{hom}^0(E,E'')
    = \mathrm{Hom}(E_{x_0},E''_{x_0})$ be the morphisms in $\mathscr
    M(L)$ determined by the short exact sequence, that is, $c_2 =
    i_{x_0}$ and $c_3 = p_{x_0}$. Then define $c_1 \in
    \mathrm{hom}^1(E'',E') =
    \bigoplus_{y}\mathrm{Hom}(E''_{y},E'_{y})$ by
    \begin{equation*}
        c_1
        = \bigoplus_y \sum_{\gamma} \pm \,
        \pi_{\gamma}' \circ q_{x_0}  \circ \pi_{\overline\gamma} \circ
        j_y,
    \end{equation*}
    where the first sum runs over all critical points $y$ of $f$ with
    Morse index 1, and the second over all gradient flow lines of $f$
    from $x_0$ to $y$, and where the sign $\pm$ is associated to
    $\gamma$ as indicated above. In fact, we view $c_1$ as living in
    $\mathrm{hom}^0(E''[-1],E')$, and $c_3$ as living in
    $\mathrm{hom}^{1}(E,E''[-1])$.

    We claim that these morphisms fit into an exact triangle
    \begin{equation*}
        E''[-1] \xrightarrow{c_1} E'\xrightarrow{c_2} E
        \xrightarrow{c_3} E''
    \end{equation*}
    in $H^0(Tw \mathscr M(L))$. In this model, the cone $C =
    \mathrm{Cone}(c_1)$ is the twisted complex
    \begin{equation*}
        C = \left(E'' \oplus E', \delta =
        \begin{pmatrix}
            0 & 0 \\
            -c_1 & 0
        \end{pmatrix}
        \right),
    \end{equation*}
    see \cite[Section 3(p)]{Seidel--Fukaya-Picard-Lefschetz-2008} (we
    suppress the shift of $c_1$). It comes together with morphisms
    $p_C = (e_{E''},0) \in \mathrm{hom}_{Tw}(C,E'')$ and $i_C =
    (0,e_{E'})^T \in \mathrm{hom}_{Tw}(E',C)$, where $e_{E''} =
    \mathrm{id}_{E_{x_0}''}$ and $e_{E'} = \mathrm{id}_{E_{x_0}'}$
    denote the chain level identity morphisms of $E'', E'$ in
    $\mathscr
    M(L)$. 

    To prove our claim, we will make use of Lemma 3.27 in
    \cite{Seidel--Fukaya-Picard-Lefschetz-2008}, which gives a
    criterion for exactness of triangles in $H^0(Tw \mathscr
    M(L))$. According to that, we have to find a cocyle $b \in
    \mathrm{hom}_{Tw}(E,C)$ such that $[b]$ is an isomorphism and
    $[\mu_{Tw}^2(p_C, b)] = [c_3]$, $[\mu_{Tw}^2(b,c_2)] = [i_C]$ in
    $H^0(Tw \mathscr M(L))$. (Again, we suppress some shifts.)
    
    We claim that $b = (b'',b')$ with $b'' = p_{x_0} \in
    \mathrm{hom}^0(E,E'')$, $b' = q_{x_0}\in \mathrm{hom}^0(E,E')$
    satisfies these requirements. The necessary computations are
    straightforward. We start by verifying that $b$ is a cocycle,
    i.e., that $\mu_{Tw}^1(b) = 0$. Unravelling the definition of
    $\mu_{Tw}^1$, cf. \cite[Section
    (3l)]{Seidel--Fukaya-Picard-Lefschetz-2008}, we obtain
    \begin{equation*}
        \mu_{Tw}^1(b) ~=~
        \begin{pmatrix}
            \mu^1 (b'')\\
            \mu^1(b') - \mu^2(c_1,b'')
        \end{pmatrix},
    \end{equation*}
    where the $\mu^d$'s are those of $\mathscr M(L)$. Now $\mu^1(b'')
    = \mu^1(p_{x_0})$ vanishes because the $p_x$ form a morphism of
    local systems. 
    To compute $\mu^2(c_1,b'')$, note that for every critical point
    $y$ there's a unique perturbed $Y$-shaped gradient tree with
    outgoing edges converging to $x_0$ and $y$ that contributes to the
    count, and that the incoming edge of this tree also converges to
    $y$ (recall that $x_0$ is the unique local minimum). In view of
    this and recalling the definition of $c_1$, we obtain
    \begin{equation*}
        \begin{aligned}
            \mu^2(c_1,b'')
            &~=~ \bigoplus_y\sum_\gamma \pm \, \pi_{\gamma}' \circ
            q_{x_0} \circ \pi_{\overline \gamma} \circ j_y \circ p_y
            \\
            &~=~ \bigoplus_y\sum_\gamma \pm\, \pi_{\gamma}' \circ
            q_{x_0} \circ \pi_{\overline \gamma} \circ
            (\mathrm{id}_{E_y} - i_y \circ q_y)\\
            &~=~\bigoplus_y \left(\sum_\gamma \pm\, \pi_{\gamma}'
                \circ q_{x_0} \circ
                \pi_{\overline \gamma} - \sum_\gamma \pm \,q_y \right)\\
            &~=~\bigoplus_y\sum_\gamma \pm\, \pi_{\gamma}' \circ q_{x_0} \circ
            \pi_{\overline \gamma}.
        \end{aligned}
    \end{equation*}
    Here we use that the $p_x$ give a morphism of local systems, i.e.,
    commute with parallel transport maps (which makes the additional
    parallel transport maps disappear that would appear in the first
    line). Moreover, we use that $\pi'_\gamma \circ q_{x_0} \circ
    \pi_{\overline \gamma} \circ i_y = \pi'_\gamma \circ q_{x_0} \circ
    i_{x_0} \circ \pi_{\overline \gamma} = \mathrm{id}_{E_y'}$, and
    that $\sum_\gamma \pm q_y = 0$ as $x_0$ is a cocyle in the usual
    Morse complex. The result of the computation is equal to
    $\mu^1(b')$, which shows that also the second component of
    $\mu^1_{Tw}(b)$ vanishes.
    
    It remains to check that $[\mu_{Tw}^2(p_C, b)] = [c_3]$,
    $[\mu_{Tw}^2(b,c_2)] = [i_C]$ in $H^0(Tw \mathscr M(L))$. After
    unravelling again definitions, the first identity follows
    immediately (on the chain level). As for the second, we obtain
    \begin{equation*}
        \mu_{Tw}^2(b,c_2) ~=~
        \begin{pmatrix}
            \mu^2(b'',c_2)\\
            \mu^2(b',c_2) - \mu^3(c_1,b'',c_2)
        \end{pmatrix}.
    \end{equation*}
    The first component is $\mu^2(b'',c_2) = p_{x_0} i_{x_0} = 0$, as
    required. As for the second component, we have $\mu^2(b',c_2) =
    q_{x_0}i_{x_0} = \mathrm{id}_{E_{x_0}'} = i_C$, and hence we are
    done if we can show that $\mu^3(c_1,b'',c_2) = 0$. Recall that
    $b'' = p_{x_0}$ and $c_2 = i_{x_0}$, and that $p_{x_0} i_{x_0} =
    0$; this together with the fact that the $i_x$ commute with
    parallel transport maps suffices to conclude that the term
    vanishes. Hence the second required identity also holds on the
    chain level.
\end{proof}

\section{Vector bundles on elliptic curves}
\label{app:vect-bundl-ellipt}
For the convenience of the reader, we collect here a couple of facts
from Atiyah's classification \cite{Atiyah--VecBunEllCurve--1957} of
vector bundles over an elliptic curve $X$ which are used in Sections
\ref{sec:homol-mirr-symm} and \ref{sec:isomorphism}.  We denote by
$\mathcal V(r,d)$ the set of isomorphism classes of vector bundles on
$X$ of rank $r$ and degree $d$.

\begin{thm}[Cf. Theorem 3 in \cite{Atiyah--VecBunEllCurve--1957}]
    \label{thm:Atiyah-Thm3}
    There exists an integer $N(r,d)$ such that for every $n \geq
    N(r,d)$, every $\mathcal E \in \mathcal V(r,d)$ fits into a short
    exact sequence
    \begin{equation*}
        0 \to \mathcal O_X^{\oplus (r-1)}(-n) \to \mathcal E \to
        \mathrm{det}\,\mathcal E((r-1)n)\to 0.
    \end{equation*}
    (Here $-(k)$ denotes tensoring by the $k^{th}$ power of the
    hyperplane bundle.)
\end{thm}

\begin{thm}[Cf. Theorem 5 in \cite{Atiyah--VecBunEllCurve--1957}]
    \label{thm:Atiyah-Thm5}
    (i) There is a unique $\mathcal F_r \in \mathcal V(r,0)$ such that
    $H^0(X;\mathcal F_r) \neq 0$, and there is a short exact sequence
    $0 \to \mathscr O_X \to \mathcal F_r \to \mathcal F_{r-1} \to 0$
    for every $r > 1$. (ii) Every $\mathcal E \in \mathcal V(r,0)$ is
    of the form $\mathcal F_r \otimes \mathcal L$ for a unique
    $\mathcal L \in \mathcal V(1,0)$.
\end{thm}

\begin{thm}[Cf. Theorems 6 and 7 in \cite{Atiyah--VecBunEllCurve--1957}]
    \label{thm:Atiyah-Thm6}
    The choice of a line bundle $\mathcal L \in \mathcal V(1,1)$
    determines natural 1-1 correspondences $\alpha_{r,d}: \mathcal
    V(h,0) \to \mathcal V(r,d)$, where $h = \mathrm{gcd}(r,d)$.  These
    are such that $\mathrm{det}\,\alpha_{r,d}(\mathcal E) =
    \mathrm{det}\,\mathcal E \otimes \mathcal L^{\otimes d}$, which
    implies that $\mathrm{det}:\mathcal V(r,d) \to \mathcal V(1,d)$
    is an $h$-1 map.

\end{thm}

Given $\mathcal L \in \mathcal V(1,1)$ and the corresponding map
$\alpha_{r,d}: \mathcal V(h,0) \to \mathcal V(r,d)$ from Theorem
\ref{thm:Atiyah-Thm6}, we write $E_{\mathcal L}(r,d) :=
\alpha_{r,d}(\mathcal F_h).$

\begin{lem}[Cf. Lemma 24 in \cite{Atiyah--VecBunEllCurve--1957}]
    \label{lem:Atiyah-lem24}
    Suppose that $\mathrm{gcd}(r,d) = 1$. Then $E_{\mathcal L}(r,d)
    \otimes \mathcal F_h \cong E_{\mathcal L}(hr,hd)$.
\end{lem}

\begin{lem}[Cf. Lemma 26 in \cite{Atiyah--VecBunEllCurve--1957}]
    \label{lem:Atiyah-lem26}
    Let $\mathcal L_1 \in \mathcal V(1,1)$. Then for every $\mathcal E
    \in \mathcal V(r,d)$ there exists some $\mathcal L_0 \in \mathcal
    V(1,0)$ such that $\mathcal E = E_{\mathcal L_1}(r,d) \otimes
    \mathcal L_0$.
\end{lem}

\bibliographystyle{alpha}
\bibliography{cobT2}

\end{document}